\documentclass[11pt,reqno]{article}

\usepackage{amsmath}
\usepackage{amssymb}
\usepackage{amsthm}
\usepackage{srcltx}
\usepackage{enumerate}
\numberwithin{equation}{section}

\newtheorem{lemma}{Lemma}[section]
\newtheorem{theorem}[lemma]{Theorem}

\newtheorem{proposition}[lemma]{Proposition}
\newtheorem{corollary}[lemma]{Corollary}

\newtheorem{remark}[lemma]{Remark}
\newtheorem{remarks}[lemma]{Remarks}

\def\Gammato{\to\!\!\!\!\!\!\!{}^\Gamma\;\;}   
\def\Gammatopi{\to\!\!\!\!\!\!\!\!\!{}^{\Gamma\widetilde\pi}\;\;}   
  
\def\2to{\to\!\!\!\!\!\!{}_{_2}\;\;}
\def\pito{\;\!\!\to\!\!\!\!\!\!\!_{_{\scriptstyle\widetilde\pi}}\;\;} 
\def\tauto{\to\!\!\!\!\!\!{}_{_\tau}\;\;}
 
\def\w2to{\rightharpoonup\!\!\!\!\!\!{}_{_2}\;\;}
\def\ws2to{\rightharpoonup\!\!\!\!\!\!{}_{_2}\!\!\!{}^{*}\;\;\;}
\def\wto{\rightharpoonup}
\def\wsto{\rightharpoonup\!\!\!\!\!\!{}^{*}\;\,}

\def\erre
{{\bf R}}

\def\enne
{{\bf N}}

\def\dps{\displaystyle}
\def\Theta{{\mit \Theta}}
\let\La=\Lambda\def\Lambda{{\mit \La}}
\let\Sg=\Sigma\def\Sigma{{\mit \Sg}}
\let\Pig=\Pi\def\Pi{{\mit \Pig}}
\let\Ps=\Psi\def\Psi{{\mit \Ps}}
\let\Xii=\Xi\def\Xi{{\mit \Xii}}
\def\graph{\mathop{\rm graph}}

\def\gr{\mathop{\rm graph}}

\title{\bf Structural stability of flows via evolutionary {\boldmath $\Gamma$}-convergence 
of weak-type}  
\author{Augusto Visintin
\thanks{Dipartimento di Matematica dell'Universit\`a degli Studi di Trento --
via Sommarive 14,  38050 Povo di Trento, Italia -- email: augusto.visintin@unitn.it }
}

\date{\today} 
\begin{document}
\maketitle 
   
\begin{abstract} 
The initial-value problem associated with multi-valued operators in Banach spaces is here reformulated as a minimization principle, extending results of Brezis-Ekeland, Nayroles and Fitzpatrick.
At the focus there is the stability of these problems w.r.t.\ perturbations not only of 
data but also of operators; this is achieved via De Giorgi's $\Gamma$-convergence 
w.r.t.\ a nonlinear topology of weak type.
A notion of evolutionary $\Gamma$-convergence of weak type is also introduced.
These results are applied to quasilinear PDEs, including doubly-nonlinear flows.  
\end{abstract}

\bigskip 
\noindent{\bf Keywords:}  
Fitzpatrick theory, 
Representative functions,
% Variational formulation,
% Null-minimization, 
Nonlinear weak convergence,
Evolutionary $\Gamma$-convergence,
Maximal monotone flows,
Doubly-nonlinear parabolic equations.
% Pseudo-monotone operators.  
 
\bigskip 
\noindent{\bf AMS Classification (2000):}  
35K60, % Boundary value problems for nonlinear parabolic PDE  
47H05, % Monotone operators (with respect to duality)  
49J40, % Variational methods including variational inequalities      
58E. % Variational problems in infinite-dimensional spaces 

\section{Introduction}\label{intro} 

\noindent 
This work deals with the behaviour of the solutions of quasilinear evolutionary PDEs, 
under perturbations not only of the data but also of the operator,
thus of the very structure of the problem.

\medskip
\noindent{\bf Structural compactness and stability.}
Stability (i.e., robustness) has an obvious applicative motivation: not only data but also 
differential operators are accessible just with some approximation;
moreover perturbations often occur, e.g.\ in the coefficients.
If properly formulated, structural stability might thus be regarded as a basic requisite for the applicative feasibility, or even the soundness, of a mathematical model. 
This point of view might also be of interest for the numerical treatment of PDEs, 
whenever differential operators are approximated by a discretization technique, 
such as finite differences, finite elements, and so on. 

In order to fix ideas, let us consider a model problem of the general form $Au\ni h$, 
$A$ being a nonlinear operator acting in a Banach space and $h$ a datum.
Given bounded families of data and of operators, we shall formulate the issue of stability 
in terms of two properties:

(i) {\it structural compactness:\/} existence of convergent sequences 
of data $\{h_n\}$ and of operators $\{A_n\}$ (in a sense to be specified), and

(ii) {\it structural stability:\/} if $A_nu_n\ni h_n$ for any $n$, 
$A_n\to A$ and $u_n\to u$, then $u$ is a solution of the asymptotic problem: $Au\ni h$.  

For (typically stationary) models represented by a minimization principle,
structural compactness and stability may adequately be formulated and addressed via 
De Giorgi's theory of $\Gamma$-convergence;
see \cite{DeFr} and e.g.\ the monographs 
\cite{At}, \cite{Bra1}, \cite{Bra2}, \cite{Bra3}, \cite{Da}.
Indirectly, this provides analogous structural properties for the associated Euler-Lagrange equations. In particular this applies to equations that are governed by a cyclically maximal monotone operator.

The issue of stability, and the techniques that we use for its analysis, 
have elements of contact with the theory of homogenization.
In that case however operators vary in a prescribed way, 
whereas here we consider arbitrary variations.
Moreover, the point of view of homogenization is quite different from the present one,
since in that case the structure of the asymptotic problem is different from that of the approximating sequence.

\medskip
\noindent{\bf Flows.} 
In the seminal paper \cite{Fi} Fitzpatrick characterized the larger class of maximal monotone operators in terms of a minimization principle.
This also provided a new interpretation of the pioneering papers of Brezis and Ekeland \cite{BrEk} and Nayroles \cite{Na}, and prompted the extension of \cite{ViAMSA}.
This method applies to flows of the form 
\begin{equation}\label{eq.intro.mm}
D_tu + \alpha(u)\ni h \qquad\hbox{ a.e.\ in time }\;(D_t := {\partial/\partial t});
\end{equation} 
here $\alpha$ is maximal monotone (e.g.\ the $p$-Laplacian) and $h$ is a known source. 
This prompts us to introduce what we shall refer to as {\it null-minimization,\/} 
see \eqref{eq.fitzp.nullmin},
and to extend the approach to stability of \cite{ViCalVar} based on $\Gamma$-convergence. 
By using a peculiar nonlinear topology, we prove the structural stability of the initial-value problem for \eqref{eq.intro.mm} and for (nonmonotone) doubly-nonlinear equations of the 
form
\begin{eqnarray} 
&D_t\partial\gamma(u) + \alpha(u)\ni h, 
\label{eq.intro.1}
\\[1mm]
&\alpha(D_tu) + \partial\gamma(u)\ni h;
\label{eq.intro.2}
\end{eqnarray}
in either case $\alpha$ is a maximal monotone operator, and $\gamma$ is
a lower semicontinuous function. 
It is known that several phenomena of mathematical physics are modeled by equations of the form
\eqref{eq.intro.mm}--\eqref{eq.intro.2}, see e.g.\ \cite{Rou2}.
 
These examples illustrate a method that may also be used to prove the structural stability of other evolutionary equations, if these are expressed in terms of maximal monotone operators.

\medskip
\noindent{\bf A nonlinear topology.} 
The theory of Fitzpatrick deals with functionals that act on the Cartesian product 
$V \!\times\! V'$ ($V$ being a Banach space and $V'$ its dual),
and are convex and lower semicontinuous.
Along with \cite{ViCalVar}, here we equip the space $V \!\times\! V'$ with
a special nonlinear topology of weak type, that we label by $\widetilde\pi$. 
This somehow exotic topology is aimed to meet two opposite exigences: 
to provide compactness of the class of functionals, and at the same time 
to yield convergence of the perturbed equation.
Compactness indeed requires a sufficiently weak topology, 
whereas one can pass to the limit in perturbed equations converge 
only if the topology is strong enough.
We use evolutionary $\Gamma$-convergence w.r.t.\ this nonlinear topology. 
 
Several papers have been dealing with stability via so-called
Mosco-convergence (see e.g.\ \cite{At}, \cite{Mo}),
and one might wonder whether this is feasible in the present set-up, too.
Although this variational convergence has a wealth of properties, 
we refrain from applying it here, since we are also concerned with compactness;
results of Mosco-compactness indeed seem to be rare,
at variance with $\Gamma$-compactness.  

\bigskip
\noindent{\bf Plan of work.}
In Section~\ref{sec.evol} we define the notion of evolutionary $\Gamma$-convergence of weak
type, and prove a theorem of $\Gamma$-compactness. 
This analysis is instrumental to Sections~\ref{sec.par}, \ref{sec.DNE} and \ref{sec.DNE'};
we isolate that result from that context in order to facilitate an independent reading. 

The remainder of this work is essentially devoted to the structural compactness and stability
of flows.  
In Sections~\ref{sec.fitzp} and \ref{sec.BEN} we revisit the Fitzpatrick theory and the related Brezis-Ekeland-Nayroles (``BEN'' for short) principle.
In Section~\ref{sec.comstab} we illustrate the structural compactness and stability of representable operators and of the associated flows in general, 
see Theorems~\ref{comsta} and \ref{teo.comp'}.

In the three other sections we apply the above results to the analysis of the structural stability 
of some initial-value problems in abstract form.
In Section~\ref{sec.par} we deal with the maximal monotone flow \eqref{eq.intro.mm}.
In Sections~\ref{sec.DNE} and \ref{sec.DNE'} we show the structural stability of initial-value 
problems for the doubly-nonlinear equations \eqref{eq.intro.1} and \eqref{eq.intro.2}. 

\bigskip
\noindent{\bf A look at the literature.} 
Fitzpatrick's Theorem~\ref{teo.Fi} was not noticed for several years, and 
was eventually rediscovered by Martinez-Legaz and Th\'era \cite{MaTh01}
and (independently) by Burachik and Svaiter \cite{BuSv02}.
Since then a rapidly expanding literature has been devoted to this theory, see e.g.\
\cite{BaBoWa}, \cite{BaWa}, \cite{Bor}, \cite{BuSv03},  
\cite{MaSv05}, \cite{MaSv08}, \cite{Pe04}, \cite{Pe04'}, besides several other articles.  
As we already pointed out, Fitzpatrick's formulation of maximal monotone flows via a null-minimization principle cast a new light upon a formulation that had been pointed out by Brezis and Ekeland \cite{BrEk} and by Nayroles \cite{Na} in 1976, prior to the Fitzpatrick theorem of 1988.
The first three authors assumed $\alpha$ cyclically maximal monotone; on the basis of the 
Fitzpatrick theory the extension to general maximal monotone operators was then rather obvious, 
see \cite{ViAMSA}. This was then further extended to nonmonotone flows in \cite{ViDNE}. 
 
Since the original formulation of \cite{BrEk} and \cite{Na} of 1976, the BEN principle was applied in several works. For the study of doubly-nonlinear evolutionary PDEs, it was used e.g.\ in \cite{RoRoSt} and \cite{St1}. 
In \cite{ViCalVar} the dependence on data and operators for the solution of quasilinear maximal monotone equations was studied, by applying $\Gamma$-convergence to the  null-minimization problem. This method was also used by this author for the homogenization of evolutionary quasilinear PDEs in other works, see e.g.\ \cite{Vi13} and references therein.
%in \cite{ViARMA}, \cite{ViJDE} and \cite{Vi13}.
In \cite{ViPse} the structural stability of pseudo-monotone equations and the corresponding doubly-nonlinear first-order flow were studied without using the Fitzpatrick theory.  

Our definition of evolutionary $\Gamma$-convergence of Section~\ref{sec.evol} may be 
reduced to $\Gamma$-conver\-gence of set-valued functionals depending on a parameter (here time). This issue was already studied in a fairly general set-up
under the denomination of $\bar\Gamma$-convergence, 
see e.g.\ \cite{Da} and references therein;
that theory includes a result of $\Gamma$-compactness (see e.g.\ Theorem~16.9 of \cite{Da}), but does not encompass the present Theorem~\ref{teo.comp}. 
A stronger notion of evolutionary $\Gamma$-convergence was addressed
in \cite{SaSe}, see also \cite{DaSa}, \cite{Mi1}, \cite{Mi2},  
where $\Gamma$-convergence was assumed for almost every instant.
Here instead we deal with weak (rather than pointwise) convergence in time, 
because of the poorness of the uniform estimates that are available.
 
The equations \eqref{eq.intro.mm}--\eqref{eq.intro.2} 
were already studied in \cite{ViDNE} and \cite{ViCalVar}, but
the present paper substantially improves the results of both works
by applying the point of view of structural compactness and stability.
\cite{ViDNE} addressed the structural stability of \eqref{eq.intro.1} and \eqref{eq.intro.2},
without considering structural compactness. 
In particular \cite{ViCalVar} concerned maximal monotone flows of the form \eqref{eq.intro.mm}, but derived the existence of a memoryless representative function only at the expense of further regularity assumptions.
This restricted the range of applicability of this technique, that indeed was just used for \eqref{eq.intro.mm}. By using a different approach, the present work removes those limitations,
and allows the extension to \eqref{eq.intro.1} and \eqref{eq.intro.2},
besides e.g.\ pseudo-monotone flows and further equations that are not addressed here.

\section{Evolutionary {\boldmath $\Gamma$}-convergence of weak type}
\label{sec.evol}  

\noindent
In this section we extend De Giorgi's basic notion of $\Gamma$-convergence to operators 
(rather than functionals) that act on time-dependent functions ranging in a Banach space $X$.  

\bigskip
\noindent{\bf Functional set-up.}
Let $X$ be a real separable and reflexive Banach space, and $p\in [1,+\infty[$. 
Let 
\begin{equation}\label{eq.evol.mu}
\begin{split} 
&\hbox{ $\mu$ be a positive and finite measure on $]0,T[$, }
\\ 
&\hbox{ which is absolutely continuous w.r.t.\ the Lebesgue measure. } 
\end{split} 
\end{equation} 
Let us set
\begin{equation}\label{eq.evol.defgamma.1}
\begin{split}
&L^p_\mu(0,T;X) = \Big\{\mu\hbox{-measurable }w: {}]0,T[{}\to X:
\!\!\int_0^T \! \|w(t)\|_X^p \, d\mu(t)< +\infty \Big\},
\\
&L^1_\mu(0,T) = \Big\{\mu\hbox{-measurable }v: {}]0,T[{}\to \erre:
\!\!\int_0^T \! |v(t)| \, d\mu(t)< +\infty \Big\},
\!\!\end{split}
\end{equation} 
and equip these spaces with the graph norm.
Examples of interest will be $\mu$ equal to the Lebesgue measure,
and $\mu$ such that $d\mu(t) = (T-t) \, dt$.

Let us equip $L^p_\mu(0,T;X)$ with a topology $\tau$ that either coincides or is finer than the weak topology.
\footnote{ We assume this in consideration of the application of 
Sections~\ref{sec.par}, \ref{sec.DNE} and \ref{sec.DNE'}.
A reader interested just in evolutionary $\Gamma$-convergence might go through this section assuming that $\mu$ is the Lebesgue measure and $\tau$ coincides with the weak topology.
}
For any operator $\psi: L^p_\mu(0,T;X)\to L^1_\mu(0,T):w\mapsto \psi_w$, let us set 
\begin{equation}\label{eq.evol.crochet}
[\psi,\xi](w) = \int_0^T \psi_w(t) \, \xi(t) \, d\mu(t) 
\qquad\forall w\in L^p_\mu(0,T;X), \forall \xi\in L^\infty(0,T). 
\end{equation}  
 
\noindent{\bf Definition of evolutionary $\Gamma$-convergence of weak type.}
Let $\{\psi_n\}$ be a bounded sequence of operators $L^p_\mu(0,T;X)\to L^1_\mu(0,T)$; 
by boundedness we mean that, for any bounded subset $A$ of $L^p_\mu(0,T;X)$, 
the set $\{\psi_{n,w}: w\in A, n\in \enne\}$ is bounded in $L^1_\mu(0,T)$.
If $\psi$ also is an operator $L^p_\mu(0,T;X)\to L^1_\mu(0,T)$, we shall say that 
\begin{equation}\label{eq.evol.defgamma.1+}
\begin{split}
&\hbox{ $\psi_n$ sequentially $\Gamma$-converges to }\psi
\\
&\hbox{ in the topology $\tau$ of $L^p_\mu(0,T;X)$ and } 
\\
&\hbox{ in the weak topology of }L^1_\mu(0,T)
\end{split}
\end{equation}
if and only if, denoting by $L^\infty_+(0,T)$ the cone of the nonnegative functions of 
$L^\infty(0,T)$,
\begin{equation}\label{eq.evol.defgamma.2}
[\psi_n, \xi] \hbox{ sequentially $\Gamma\tau$-converges to $[\psi,\xi]$ in }
L^p_\mu(0,T;X), \; \forall \xi\in L^\infty_+(0,T).
\end{equation}  
We shall say that a sequence {\it $\Gamma\tau$-converges\/} if it $\Gamma$-converges w.r.t.\ $\tau$, that a functional is {\it $\tau$-lower semicontinuous\/} if it is lower semicontinuous w.r.t.\ $\tau$, and so on.

By the classical definition of sequential $\Gamma$-convergence, 
\footnote{ \eqref{eq.evol.defgamma.3} is often referred to as the {\it inferior-limit condition;\/} 
the sequence $\{w_n\}$ occurring in \eqref{eq.evol.defgamma.4} is called a 
{\it recovery sequence.\/} 
}
\eqref{eq.evol.defgamma.2} means that for any $\xi\in L^\infty_+(0,T)$ 
\begin{equation}\label{eq.evol.defgamma.3}
\begin{split}
&\forall w\in L^p_\mu(0,T;X), \forall \hbox{ sequence $\{w_n\}$ in }L^p_\mu(0,T;X), 
\\
&\hbox{if \ $w_n\tauto w$ in $L^p_\mu(0,T;X)$ \ then \ }
\liminf_{n\to+ \infty} \; [\psi_n,\xi](w_n) \ge [\psi,\xi](w), 
\end{split}
\end{equation}
  \vskip-0.4truecm
\begin{equation}\label{eq.evol.defgamma.4}
\begin{split}
&\forall w\in L^p_\mu(0,T;X), \exists\hbox{ sequence $\{w_n\}$ 
of $L^p_\mu(0,T;X)$ such that } 
\\
&w_n\tauto w \hbox{ in $L^p_\mu(0,T;X)$ \ and \ }
\lim_{n\to +\infty} \; [\psi_n,\xi](w_n) = [\psi,\xi](w). 
\end{split}
\end{equation}

By the properties of ordinary $\Gamma$-convergence, \eqref{eq.evol.defgamma.2} entails that
\begin{equation}\label{eq.evol.defgamma.5}
[\psi, \xi] \hbox{ is sequentially $\tau$-lower semicontinuous in }L^p_\mu(0,T;X), \; 
\forall \xi\in L^\infty_+(0,T).
\end{equation}
 
\begin{remarks}\rm
(i) This definition of evolutionary $\Gamma$-convergence is not equivalent either to that of \cite{SaSe}, or to that of \cite{DaSa}, \cite{Mi1}, \cite{Mi2}.
In those works $\Gamma$-convergence is actually assumed for almost any 
$t\in {}]0,T[$, whereas here it is just weak in $L^1_\mu(0,T)$. 

(ii) The present definition of parameter-dependent $\Gamma$-convergence 
is based on testing the sequence $\{\psi_n\}$ on functions of time. 
Ahead in this section we shall see that this notion might equivalently be reformulated in terms of set-valued functions. 

Apart from this minor difference, this set-up fits the general framework of 
$\bar\Gamma$-convergence, defined in Chap.~16 of \cite{Da}, see also references therein.
However, here we deal with functions of the parameter that are integrable w.r.t.\ the Lebesgue measure, rather than w.r.t.\ more general measures as in \cite{Da}.

(iii) The deep analysis of Chaps.~16-20 of \cite{Da} does not seem to encompass 
the next theorem, that rests upon a result of \cite{Hi}. 
\end{remarks}
  
Although we considered generic operators $\psi: L^p(0,T;X)\to L^1(0,T)$,
our main concern is for memoryless operators of the form  
\begin{equation}\label{eq.evol.super}
\begin{split}
&\psi_w(t) = \varphi(t,w(t))  
\qquad\forall w\in L^p_\mu(0,T;X),\hbox{ for a.e.\ }t\in {}]0,T[,
\\
&\varphi: {}]0,T[{} \times X\to \erre^+ \hbox{ \ being a {\it normal function.\/} }
\end{split}
\end{equation} 
By this we mean that $\varphi$ is globally measurable and $\varphi(t,\cdot)$ is lower semicontinuous for a.e.\ $t\in {}]0,T[$. 
  
\begin{theorem} \label{teo.comp}
Let $X$ be a real separable and reflexive Banach space, $p\in [1,+\infty[$,
and $\{\varphi_n\}$ be a sequence of normal functions 
${}]0,T[{} \times X\to \erre^+$.
Assume that this sequence is equi-coercive and equi-bounded, in the sense that  
\begin{equation}\label{eq.evol.equibc} 
\begin{split}
&\exists C_1,C_2,C_3 >0: \forall n,\hbox{for a.e.\ }t\in {}]0,T[, \forall w\in X, 
\\
&C_1 \|w\|_X^p\le \varphi_n(t,w) \le C_2\|w\|_X^p +C_3, 
\end{split} 
\end{equation} 
and that
\begin{equation}\label{eq.evol.nul}
\varphi_n(t,0) =0 \qquad\hbox{ for a.e.\ }t\in {}]0,T[,\forall n.
\end{equation}  
\indent
Let $\mu$ fulfill \eqref{eq.evol.mu}. Let
$\tau$ be a topology on $L^p_\mu(0,T;X)$ that either coincides or is finer 
than the weak topology, and such that 
\begin{equation}\label{eq.evol.gamcom} 
\begin{split}
&\hbox{ for any sequence $\{F_n\}$ of functionals }
L^p_\mu(0,T;X)\to \erre^+\cup \{+\infty\},
\\
&\hbox{ if \ }\sup_{n\in\enne} \big\{\|w\|_{L^p_\mu(0,T;X)}: 
w\in L^p_\mu(0,T;X), F_n(w) \le C \big\}<+\infty,
\\
&\hbox{ then \ $\{F_n\}$ has a sequentially $\Gamma\tau$-convergent subsequence. }  
\end{split}
\end{equation}
\indent
Then there exists a normal function $\varphi: {}]0,T[{} \times X\to \erre^+$ such that
$\varphi(\cdot,0) =0$ a.e.\ in $]0,T[$, and such that, 
defining the operators $\psi,\psi_n: L^2_\mu(0,T;X)\to L^1_\mu(0,T)$ for any $n$ 
as in \eqref{eq.evol.super}, possibly extracting a subsequence,
\begin{equation}\label{eq.evol.tesi} 
\begin{split}
&\hbox{ $\psi_n$ sequentially $\Gamma$-converges to }\psi
\\
&\hbox{ in the topology $\tau$ of $L^p_\mu(0,T;X)$ and } 
\\
&\hbox{ in the weak topology of $L^1_\mu(0,T)$ (cf.\ \eqref{eq.evol.defgamma.2}). }
\end{split}
\end{equation}
\indent
Moreover, if $\varphi_n$ does not explicitly depend on $t$ for any $n$, 
then the same holds for $\varphi$.
\end{theorem}  

\noindent{\bf Proof.\/}  
For the reader's convenience, we split this argument into several steps.

(i) First we show that, denoting by $C^0_+([0,T])$ the cone of the nonnegative functions of 
$C^0([0,T])$, 
\begin{equation}\label{eq.evol.conv}
\begin{split}
&\forall \xi\in C^0_+([0,T]), \exists g_\xi: L^p_\mu(0,T;X)\to \erre^+ 
\hbox{ such that } 
\\
&[\psi_n,\xi] \hbox{ $\Gamma\tau$-converges to }g_\xi.
\end{split}
\end{equation}

The separable Banach space $C^0([0,T])$ has a countable dense subset $M$,
e.g., the polynomials with rational coefficients. Let us
define the cone $M_+$ of the nonnegative elements of $M$.
For any $\xi\in M_+$, by \eqref{eq.evol.equibc} a suitable subsequence $\{[\psi_{n'},\xi]\}$ 
weakly $\Gamma$-converges to a function $g_\xi: L^p_\mu(0,T;X)\to \erre$, and
\begin{equation}\label{eq.evol.a}
\begin{split}
&C_1 \int_0^T \|w(t)\|_X^p \, \xi(t) \, d\mu(t) \le g_\xi(w) 
\\
&\le C_2\int_0^T \|w(t)\|_X^p \, \xi(t) \, d\mu(t) + C_3 \int_0^T \xi(t) \, d\mu(t)
\qquad\forall w\in L^p_\mu(0,T;X). 
\end{split}
\end{equation}  

A priori the selected subsequence $\{[\psi_{n'},\xi]\}$ might depend on $\xi$. 
However, because of the countability of $M_+$, via a diagonalization procedure 
one may select a subsequence that is independent of $\xi\in M_+$.
(Henceforth we shall write $\psi_n$ in place of $\psi_{n'}$, dropping the prime.) 
For that subsequence thus
\begin{equation}\label{eq.evol.a=}
\begin{split}
&\forall \xi\in M_+, \exists g_\xi: L^p_\mu(0,T;X)\to \erre^+ \hbox{ such that } 
\\
&[\psi_n,\xi]\hbox{ \ $\Gamma\tau$-converges to \ $g_\xi$ in }L^p_\mu(0,T;X),
\end{split}
\end{equation}
that is, for any $\xi\in M_+$,
\begin{equation}
\begin{split}\label{eq.evol.b=}
&\forall w\in L^p_\mu(0,T;X), \forall \hbox{ sequence $\{w_n\}$ in }L^p_\mu(0,T;X), 
\\
&\hbox{if \ $w_n \tauto w$ in $L^p_\mu(0,T;X)$ \ then \ }
\liminf_{n\to+ \infty} \; [\psi_n,\xi](w_n) \ge g_\xi(w), 
\end{split}
\end{equation}
       \vskip-0.4truecm 
\begin{equation}
\begin{split}\label{eq.evol.c=}
&\forall w\in L^p_\mu(0,T;X), \exists\hbox{ sequence $\{w_n\}$ of $L^p_\mu(0,T;X)$ 
such that } 
\\
&w_n \tauto w \hbox{ in $L^p_\mu(0,T;X)$ \ and \ }
\lim_{n\to + \infty} \; [\psi_n,\xi](w_n) = g_\xi(w). 
\end{split}
\end{equation} 
(The recovery sequence $\{w_n\}$ in \eqref{eq.evol.c=} may depend on $\xi$.)

As any $\xi\in C^0_+([0,T])$ is the uniform limit of some sequence $\{\xi_m\}$ in $M$,
for any bounded sequence $\{w_n\}$ in $L^p_\mu(0,T;X)$
\begin{equation}
\begin{split}
&\sup_n \big| [\psi_n,\xi](w_n) - [\psi_n,\xi_m](w_n) \big| 
\overset{\eqref{eq.evol.crochet}}{=}
\sup_n \Big|\int_0^T \psi_{n,w_n}(t) \, [\xi(t) - \xi_m(t)] \, d\mu(t) \Big| 
\\
&\overset{\eqref{eq.evol.super}}{\le}
\| \xi - \xi_m \|_{C^0([0,T])} \sup_n \int_0^T \varphi_n(t,w) \, d\mu(t) 
\\
&\overset{\eqref{eq.evol.equibc}}{\le}\| \xi - \xi_m \|_{C^0([0,T])}
\sup_n \bigg\{ C_2 \!\! \int_0^T \|w_n(t)\|_X^p \, d\mu(t) 
+ C_3 \mu(]0,T[) \bigg\}
\qquad\forall m.
\end{split}
\end{equation}
By the density of $M_+$ in $C^0_+([0,T])$, the convergence \eqref{eq.evol.a=} 
then holds for any $\xi\in C^0_+([0,T])$.

\medskip
(ii) Next we extend the $\Gamma\tau$-convergence \eqref{eq.evol.conv} 
to any $\xi\in L^\infty_+(0,T)$.

By the classical Lusin theorem, for any $\xi\in L^\infty_+(0,T)$ there exists a sequence 
$\{\xi_m\}$ in $C^0_+([0,T])$ such that  
\begin{equation}\label{eq.evol.Lusin}
\begin{split}
&\|\xi_m\|_{C^0([0,T])} \le \|\xi\|_{L^\infty(0,T)} \quad\forall m, \hbox{ and}
\\[1mm]
&\hbox{setting }A_m =\{t\in [0,T]: \xi_m(t) \not=  \xi(t)\}, \quad \mu(A_m)\to 0.
\end{split}
\end{equation}  
Hence 
\begin{equation}
\begin{split}
&\sup_n \big| [\psi_n,\xi](w_n) - [\psi_n,\xi_m](w_n) \big| 
\overset{\eqref{eq.evol.crochet}}{=} 
\sup_n \Big|\int_0^T \psi_{n,w_n}(t) \, [\xi(t) - \xi_m(t)] \, d\mu(t) \Big| 
\\ 
&\overset{\eqref{eq.evol.super},\eqref{eq.evol.equibc}}{\le} \|\xi - \xi_m\|_{L^\infty(0,T)} 
\bigg\{ C_2 \!\! \int_{A_m} \|w_n(t)\|_X^p \, d\mu(t) + C_3 \mu(A_m) \bigg\}
\qquad\forall m.
\end{split}
\end{equation}
As $w_n\tauto w$ in $L^p_\mu(0,T;X)$ (see \eqref{eq.evol.c=}), 
the sequence $\{\|w_n(\cdot)\|_X^p\}$ is equi-integrable. 
By this property and by \eqref{eq.evol.Lusin}, the latter expression vanishes as $m\to \infty$.
\eqref{eq.evol.conv} is thus extended to any $\xi\in L^\infty_+(0,T)$. 

\medskip
(iii) Next we prove that  
\begin{equation}\label{eq.evol.due}
\begin{split}
&\exists \psi: L^p_\mu(0,T;X)\to L^1_\mu(0,T) \hbox{ such that } 
\\
&g_\xi(w) = [\psi_w,\xi] \qquad\forall w\in L^p_\mu(0,T;X),\forall \xi\in L^\infty_+(0,T).
\end{split}
\end{equation} 
(Some care is needed, since $L^\infty_+(0,T)$ is not a linear space.)
Let us fix any $\xi\in L^\infty_+([0,T])$, any $w\in L^p_\mu(0,T;X)$, and any sequence
$\{w_n\}$ as in \eqref{eq.evol.c=}. 
By the boundedness of $\{w_n\}$ and by \eqref{eq.evol.equibc}, 
the sequence $\{\psi_{n,w_n}\}$ = $\{\varphi_n(\cdot,w_n)\}$ 
is bounded in $L^1_\mu(0,T)$ and is equi-integrable. 
There exists then a function $\gamma \in L^1_\mu(0,T)$ such that, 
possibly extracting a subsequence,
$\psi_{n,w_n}\wto \gamma$ weakly in $L^1_\mu(0,T)$. 
\footnote{ We denote the strong and weak convergence respectively by $\to$ and $\wto$. 
} 
Thus
\begin{equation}\label{eq.evol.o}
[\psi_n,\xi](w_n) =\int_0^T \psi_{n,w_n}(t) \, \xi(t) \, d\mu(t) 
\to \int_0^T \gamma(t) \, \xi(t) \, d\mu(t)
\qquad\forall \xi\in L^\infty(0,T).
\end{equation}  
Thus by \eqref{eq.evol.c=}
\begin{equation}\label{eq.evol.q}
g_\xi(w) = \int_0^T \gamma(t) \, \xi(t) \, d\mu(t)
\qquad\forall \xi\in L^\infty_+(0,T).
\end{equation}

Therefore $\gamma$ is determined by $w\in L^p_\mu(0,T;X)$ 
(and of course by the sequence  $\{\psi_n\}$), 
but is independent of the specific sequence $\{w_n\}$ that fulfills \eqref{eq.evol.c=}.  
This defines an operator 
\begin{equation}\label{eq.evol.f}
\psi: L^p_\mu(0,T;X)\to L^1_\mu(0,T): w \mapsto \psi_w = \gamma.
\end{equation} 
The equality \eqref{eq.evol.q} thus reads
\begin{equation}\label{eq.evol.r}
g_\xi(w) =\int_0^T \psi_w(t) \, \xi(t) \, d\mu(t) 
\qquad\forall \xi\in L^\infty_+(0,T), \forall w\in L^p_\mu(0,T;X).
\end{equation} 
Recalling the definition \eqref{eq.evol.crochet}, 
we see that this completes the proof of \eqref{eq.evol.due}.  

\medskip
(iv) Finally, we show that there exists a normal function 
$\varphi: {}]0,T[{} \times X\to \erre^+$ such that the operator
$\psi$ that we just defined in \eqref{eq.evol.f} is as in \eqref{eq.evol.super}. 

By \eqref{eq.evol.nul}, for any $n$ the functional
\begin{equation}
\begin{split}
&\Phi_n: L^p_\mu(0,T;X)\to \erre^+: 
\\
&w\mapsto \int_0^T \psi_{w,n}(t) \, d\mu(t) =\int_0^T \varphi_n(t,w_n(t)) \, d\mu(t)
\end{split}
\end{equation}
is additive in the sense of \eqref{eq.evol.add} below. 
This property then also holds for the limit functional
\begin{equation}
\Phi: L^p_\mu(0,T;X)\to \erre^+: w\mapsto \int_0^T \psi_w(t) \, d\mu(t).
\end{equation}
By selecting $\xi\equiv 1$ in \eqref{eq.evol.defgamma.5}, we get that $\Phi$ is 
lower semicontinuous. 
By Lemma~\ref{lemma.Hiai} below then there exists a normal function $\varphi$ as we just specified.
\hfill$\Box$  

\bigskip
As a preparation for the next lemma,
let us say that a functional $\Phi: L^p_\mu(0,T;X)\to \erre$ is invariant by translations if
$\Phi(\widetilde w(\cdot +s)) = \Phi(w)$
whenever $w$ and $s$ fulfill the following conditions: 
$s>0$, $w\in L^p_\mu(0,T;X)$ and 
(denoting by $\widetilde w$ the function obtained by extending $w$ to $\erre$ 
with null value) $\widetilde w(t +s)=0$ for a.e.\ $t\in \erre\setminus {} ]0,T[$. 

\begin{lemma} [\cite{Hi}] \label{lemma.Hiai}
Let $X$ be a real separable Banach space and $p\in [1,+\infty[$.
Let a functional $\Phi: L^p_\mu(0,T;X)\to \erre\cup \{+\infty\}$ be lower semicontinuous, 
proper (i.e., $\Phi \not\equiv +\infty$) and {\rm additive,\/} in the sense that  
\begin{equation}\label{eq.evol.add}
\begin{split} 
&\forall w_1,w_2\in L^p_\mu(0,T;X), 
\\
&\mu\big( \{t\in {}]0,T[{}: w_1(t)w_2(t) \not= 0\} \big) =0
\quad\Rightarrow\quad \Phi(w_1 +w_2) = \Phi(w_1) + \Phi(w_2).
\end{split}
\end{equation}
\indent
Then there exists a normal function $\varphi: {}]0,T[{} \times X\to \erre\cup\{+\infty\}$ 
such that
\begin{eqnarray}
\begin{split}\label{aaa}
&\varphi(t,\cdot) \not\equiv +\infty \quad\forall t\in ]0,T[,
\\
&\varphi(\cdot,0) =0 \quad\hbox{ a.e.\ in }{}]0,T[,
\\ 
&\Phi(w) = \int_0^T \varphi(t, w(t)) \, d\mu(t)
\qquad\forall w\in L^p_\mu(0,T;X).
\end{split} 
\end{eqnarray}  
Moreover,
\begin{equation}\label{eq.evol.add}
\hbox{ if $\Phi$ is convex then $\varphi(t,\cdot)$ is also convex 
for a.e.\ $t\in{} ]0,T[$, }
\end{equation} 
and the function $\varphi$ is unique, up to sets of the form $N\times X$ with $\mu(N)=0$.
\\
\indent
Finally, if the functional $\Phi$ is invariant by translations, 
then $\varphi$ does not explicitly depend on $t$.
\end{lemma}  

\begin{remarks}\rm \label{rem.amend}
Theorem~\ref{teo.comp} may be applied if $\tau$ coincides with the weak topology 
or with the topology $\widetilde\pi$ of $L^p_\mu(0,T;X)$, 
that we shall introduce in Section~\ref{sec.comstab} ahead.
In either case the hypothesis \eqref{eq.evol.gamcom} is fulfilled, see Corollary~8.18 of \cite{Da} for the weak topology, and Theorem~4.4 of \cite{ViCalVar} for the topology 
$\widetilde\pi$. 

(ii) Denoting by ${\cal L}(0,T)$ the $\sigma$-algebra of the 
Lebesgue-measurable subsets of $]0,T[$, 
\eqref{eq.evol.defgamma.2} is equivalent to
\begin{equation}\label{eq.evol.defgamma.2'}
\int_A \psi_{n,w}(t) \, d\mu(t) \to \int_A \psi_w(t) \, d\mu(t) 
\qquad\forall A\in {\cal L}(0,T).
\end{equation} 
As the sets of ${\cal L}(0,T)$ are in one-to-one correspondence 
with the characteristic functions of $L^\infty(0,T)$,
this equivalence may be proved by mimicking the argument based on the Lusin theorem,
that we used in the proof above. 

(iii) Lemma~\ref{lemma.Hiai} may be compared e.g.\ with Section~2.4 of \cite{But},
which deals with a finite-dimensional space $X$. 
\end{remarks}

\section{Brief review of Fitzpatrick's theory}
\label{sec.fitzp}

\noindent 
In this section we review a celebrated result of Fitzpatrick, by which any maximal monotone 
operator can be formulated as a minimization principle, 
and briefly illustrate some results of that theory. 
 
\bigskip
\noindent{\bf The Fitzpatrick Theorem.} 
Let $V$ be a real Banach space, and denote by $\langle \cdot,\cdot \rangle$ 
the duality pairing between $V'$ and $V$. 
Let $\alpha: V \to {\cal P}(V')$ be a (possibly multi-valued) measurable operator,
i.e., such that
\begin{equation}
\begin{split}
&g^{-1}(A) := \big\{v\in V: g(v)\cap A\not= \emptyset \big\}
\\[1mm]
&\hbox{is measurable, for any open subset $A$ of $V'$. }
\end{split}
\end{equation}
For instance, this condition is fulfilled if $\alpha$ is maximal monotone.
Let us assume that $\alpha$ is proper, i.e., $\alpha(V) \not= \emptyset$. 

In \cite{Fi} Fitzpatrick defined what is now called the {\it Fitzpatrick function:\/}
\begin{equation}\label{eq.fitzp.Fitzfun}
\begin{split}
f_\alpha(v,v^*) :=
& \; \langle v^*,v\rangle + \sup \big\{ \langle v^*- v^*_0,v_0 - v\rangle:
v^*_0\in \alpha(v_0) \big\} 
\\
=& \sup \big\{ \langle v^*, v_0\rangle - \langle v^*_0, v_0- v\rangle:
v^*_0\in \alpha(v_0) \big\} \quad\forall (v,v^*)\in V \!\times\! V'.
\end{split}  
\end{equation} 
This function is convex and lower semicontinuous, as it is
the supremum of a family of affine and continuous functions.

\begin{theorem}[\cite{Fi}]\label{teo.Fi}
If $\alpha: V \to {\cal P}(V')$ is maximal monotone then
\begin{eqnarray}
&f_\alpha(v,v^*) \ge \langle v^*,v\rangle
\qquad\forall (v,v^*)\in V \!\times\! V', 
\label{eq.fitzp.Fitztheo1}
\\
&f_\alpha(v,v^*) = \langle v^*,v\rangle
\quad\Leftrightarrow\quad v^*\in \alpha(v). 
\label{eq.fitzp.Fitztheo2}
\end{eqnarray} 
\end{theorem}

This theorem generalizes the following classical Fenchel result of convex analysis. 
Let $\varphi: V\to \erre\cup \{+\infty\}$ be a convex and lower semicontinuous 
proper function, and denote by
$\varphi^*: V'\to \erre\cup \{+\infty\}$ and $\partial\varphi: V\to {\cal P}(V')$ 
respectively the convex conjugate function and the subdifferential of $\varphi$. 
% (see e.g.\ \cite{EkTe},\cite{Fe},\cite{HiLe},\cite{Ro1}). 
Then
\begin{equation}\label{eq.fitzp.Fen}
\begin{split} 
&\varphi(v)+ \varphi^*(v^*) \ge \langle v^*,v\rangle
\qquad\forall (v,v^*)\in V \!\times\! V', 
\\
&\varphi(v)+\varphi^*(v^*)= \langle v^*,v\rangle
\;\;\Leftrightarrow\;\;
v^*\in\partial \varphi(v). 
\end{split} 
\end{equation} 
We shall refer to \eqref{eq.fitzp.Fitztheo1}, \eqref{eq.fitzp.Fitztheo2} 
(\eqref{eq.fitzp.Fen}, resp.) as the 
{\it Fitzpatrick system\/} (the {\it Fenchel system,\/} resp.), and to the mapping
$(v,v^*)\mapsto \varphi(v)+\varphi^*(v^*)$ as a {\it Fenchel function.\/}
 
\bigskip 
\noindent{\bf Representative functions.}
Extending the theory of Fitzpatrick, nowadays one says that a function $f$ (variationally) 
{\it represents\/} a proper measurable operator $\alpha:V \to {\cal P}(V')$ whenever 
\begin{equation}\label{eq.fitzp.convrep}
\begin{split} 
&f:V \!\times\! V' \to \erre\cup \{+\infty\} 
\hbox{ \ is convex and lower semicontinuous, }
\\
&f(v,v^*) \ge \langle v^*,v\rangle
\qquad\forall (v,v^*)\in V \!\times\! V',
\\
&f(v,v^*) = \langle v^*,v\rangle
\quad\Leftrightarrow\quad v^*\in \alpha(v).
\end{split}
\end{equation}  
One accordingly says that $\alpha$ is {\it representable,\/} and 
that $f$ is a {\it representative function.\/}
We shall denote by ${\cal F}(V)$ the class of the functions that 
fulfill the first two of these properties.
Representable operators are monotone, see \cite{Fi}, 
but need not be either cyclically monotone or maximal monotone.

For instance, let us define the convex and continuous mapping
\begin{equation}\label{eq.fitzp.counterex}
f_b: \erre^2\to \erre: (x,y)\mapsto b(x^2+y^2) 
\qquad\forall b>0.
\end{equation}
For $b=1/2$ this function represents the identity operator $\erre\to \erre$, 
whereas for any other $b>0$ it represents the operator 
$\alpha(0) =\{0\}$, $\alpha(x) = \emptyset$ for any $x\not=0$.

For any convex and lower semicontinuous function 
$\varphi: V\to \erre\cup \{+\infty\}$, the Fenchel function
$f: (v,v^*)\mapsto \varphi(v)+ \varphi^*(v^*)$ represents the operator 
$\partial\varphi$ because of \eqref{eq.fitzp.Fen}.
(Incidentally note that this coincides with the Fitzpatrick function 
$f_{\partial\varphi}$ only exceptionally).
Other examples of interest for the theory of PDEs were provided e.g.\ in \cite{ViDNE} and 
\cite{ViCalVar}; they include for instance a representative function for quasilinear elliptic operators of the form
\begin{equation}
H^1_0(\Omega)\to {\cal P}(H^{-1}(\Omega)): 
v\mapsto - \nabla \cdot \vec\phi(\nabla v)
\qquad(\Omega \hbox{ being a domain of }\erre^N)
\end{equation}  
for a maximal monotone mapping $\vec\phi: \erre^N\to {\cal P}(\erre^N)$.
 
Henceforth we shall assume that the real Banach space $V$ is also reflexive ---
an assumption that provides a number of results.
Besides the duality between $V$ an $V'$, we shall consider 
the duality between the product space $V \!\times\! V'$ and its dual $V' \!\times\! V$,
and the corresponding convex conjugation. More specifically,
for any function $g:V \!\times\! V' \to \erre\cup \{+\infty\}$, we define
$g^*:V' \!\times\! V \to \erre\cup \{+\infty\}$ by setting
\begin{equation}
\begin{split}
g^*(w^*,w) :=
\sup \big\{ \langle w^*,v\rangle + \langle v^*,w\rangle - g(v,v^*) : 
(v,v^*)\in V \!\times\! V'\big\} &
\\
\qquad\forall (w^*,w)\in V' \!\times\! V. &
\end{split}
\end{equation} 

Let us denote by $I_\alpha: V \!\times\! V' \to \erre\cup \{+\infty\}$ the indicator function
of the graph of any operator $\alpha:V \to {\cal P}(V')$,
i.e., $I_\alpha(v,v^*) =0$ if $v^*\in \alpha(v)$, $I_\alpha(v,v^*) =+\infty$ otherwise.
Let us also denote by $\pi$ the mapping associated to the duality pairing:
\begin{equation}\label{eq.fitzp.pi}
\pi: V \!\times\! V'\to \erre: (v,v^*)\mapsto \langle v^*,v\rangle,
\end{equation}
and by ${\cal I}$ the permutation operator
$V \!\times\! V'\to V' \!\times\! V: (v,v^*)\mapsto (v^*,v)$.
The definition \eqref{eq.fitzp.Fitzfun} thus also reads 
$f_\alpha= (\pi+ I_\alpha)^* \circ{\cal I}$, that is, 
\begin{equation}\label{eq.fitzp.Fitzfun1}
f_\alpha(v,v^*) = (\pi+ I_\alpha)^*(v^*,v)  
\qquad\forall (v,v^*)\in V \!\times\! V'.
\end{equation} 
 
Next we state some relevant results about representative functions. 
  
\begin{theorem} [\cite{BuSv03}, \cite{Sv03}]\label{teo.maxm}
Let $V$ be a reflexive Banach space. 
A function $g\in {\cal F}(V)$ represents a maximal monotone operator
$\alpha:V \to {\cal P}(V')$ if and only if \ $g^*\in {\cal F}(V')$. 
If this holds, then $g^*$ represents the inverse operator $\alpha^{-1}:V' \to {\cal P}(V)$,
that is, $g^*\circ{\cal I}$ ($\in {\cal F}(V)$) represents $\alpha$.
\end{theorem} 

Because of this theorem, if $\alpha$ is maximal monotone, then 
\[
(f_\alpha\circ{\cal I})^* = (\pi+ I_\alpha)^{**} \text{ \ \ also represents }\alpha.
\]

\begin{theorem} [\cite{BuSv02}, \cite{Fi}, \cite{MaSv08}, \cite{Pe04}] \label{teo.rep}
Let $V$ be a reflexive Banach space. 
Let $\alpha:V \to {\cal P}(V')$ be a maximal monotone operator,
$f_\alpha$ be its Fitzpatrick function, and
$g: V \!\times\! V'\to \erre\cup \{+\infty\}$ be a convex 
and lower semicontinuous function. Then
\begin{equation}\label{eq.fitzp.rep}
\begin{split}
&g\in{\cal F}(V), \quad g \hbox{ represents }\alpha \quad\Leftrightarrow\quad 
\\[1mm]
&f_\alpha \le g\le f_\alpha^*\circ{\cal I} 
\;\;\hbox{ pointwise in }V \!\times\! V'. 
\end{split}
\end{equation}
\end{theorem}   
 
\begin{theorem} [\cite{Sv03}]
Let $V$ be a reflexive Banach space. 
Any maximal monotone operator $\alpha:V \to {\cal P}(V')$ may be represented by a function 
$g\in{\cal F}(V)$ such that $g^* = g\circ{\cal I}^{-1}$. 
\end{theorem} 

An explicit examples of such a {\it self-dual\/} representative was first provided in \cite{BaWa}.
 
\bigskip 
\noindent{\bf Two minimization principles.} 
Let a function $g\in {\cal  F}(V)$ represent a proper operator
$\alpha: V\to {\cal P}(V')$, and let us define the function 
\begin{equation}\label{eq.fitzp.nullmindef}
J(v,v^*) := f(v,v^*) - \langle v^*,v\rangle 
\qquad\forall (v,v^*)\in V \!\times\! V'.
\end{equation}
Minimizing $J$ w.r.t.\ both variables $v,v^*$, the system \eqref{eq.fitzp.convrep} yields 
\begin{equation}\label{eq.fitzp.equi}
J(v,v^*) =\inf J 
\quad\Leftrightarrow\quad 
J(v,v^*) =0
\quad\Leftrightarrow\quad 
v^*\in \alpha(v).
\end{equation}  

If instead we minimize $J(v,v^*)$ just w.r.t.\ $v$ for any fixed $v^*\in V'$, then
\begin{equation}\label{eq.fitzp.nullmin}
J(v,v^*) =\inf J(\cdot,v^*) =0
\quad\Leftrightarrow\quad v^*\in \alpha(v).
\end{equation}  
Here in the implication ``$\Rightarrow$'' it is necessary to prescribe the minimum value to be zero, since it is not guaranteed that $v^*\in \alpha(V)$. 

We shall speak of {\it null-minimization\/} for a minimization problem in which 
the minimum value is prescribed to be zero as in \eqref{eq.fitzp.nullmin}, 
of {\it ordinary minimization\/} otherwise as in \eqref{eq.fitzp.equi}.

The relevance of prescribing the minimum value may conveniently be understood even 
on a linear first-order ordinary differential equation, along the lines of \cite{Vi14}.
The corresponding Euler-Lagrange equation is of second order in time, and thus 
it is not determined just by an initial condition.
Prescribing the minimum of the functional to be zero may provide the required further side-condition. In any case, existence must be proved even for a null-minimizer.

\section{Extended Brezis-Ekeland-Nayroles principle}
\label{sec.BEN}

\noindent 
In this section we outline an extension of what we shall refer to as the {\it Brezis-Ekeland-Nayroles (BEN) principle.\/}
This will be the basis for a further formulation of the next section. 

\begin{lemma} \label{lem1} 
Let $X$ be a Hausdorff topological space,
$Y$ be the dual of a Banach space equipped with the weak star topology, and
$\varphi: X \!\times\! Y\to \erre\cup \{+\infty\}$ be lower semicontinuous. If
\begin{equation}\label{eq.fitzp.coee}
\inf_{x\in X} \; \varphi(x,y) \to +\infty
\quad\hbox{ as }\|y\|_Y \to +\infty, 
\end{equation} 
then the function 
$X\to \cup \{+\infty\}: x\mapsto \inf_{y\in Y} \varphi(x,y)$ 
is lower semicontinuous.  
\end{lemma}  

This result stems from the classical Banach-Alaoglu theorem and from the next lemma,
since \eqref{eq.fitzp.coee} allows one to confine $y$ to a bounded 
(namely, relatively weakly compact) subset of $Y$.

\begin{lemma} [\cite{AuEk} p.\ 120 e.g.] \label{lem2} 
Let $X$ be a Hausdorff topological space,
$Y$ be a compact space, and
$\varphi: X \!\times\! Y\to \erre\cup \{+\infty\}$ be a lower semicontinuous function.
The function $X\to \erre: x\mapsto \inf_{y\in Y} \varphi(x,y)$ 
is then lower semicontinuous. 
\end{lemma} 
 
\begin{proposition} \label{sum}
Let $V$ be a real reflexive Banach space, and assume that
\begin{eqnarray}
&\alpha_i :V \to {\cal P}(V') \hbox{ is proper and representable }\; (i=1,2),
\\
&\exists \bar v\in V: \alpha_1(\bar v)\not=\emptyset \hbox{ and } 
\alpha_2(\bar v) \not=\emptyset,
\\
&g_i\in {\cal F}(V), \hbox{ and $g_i$ represents }\alpha_i \;(i=1,2),
\\
&\displaystyle \lim_{\|z^*\|_{V'} \to +\infty} \;
\inf_{(v,v^*)\in V \!\times\! V'} \big\{g_1(v,v^*-z^*) + g_2(v,z^*)\big\} = +\infty.
\label{eq.fitzp.coo}
\end{eqnarray}
\indent
Then: 
(i) The {\rm partial inf-convolution\/} 
\begin{equation}\label{eq.fitzp.infcon}
g_1 \!\oplus\! g_2 (v,v^*) := 
\inf_{z^*\in V'} \big\{g_1(v,v^*-z^*) + g_2(v,z^*)\big\}
\qquad\forall (v,v^*)\in V \!\times\! V'
\end{equation}
is an element of ${\cal F}(V)$, and represents the operator $\alpha_1 + \alpha_2$.
\\
\indent
(ii) If moreover $\alpha_2:V\to V'$ is linear and continuous, then 
\begin{equation}\label{eq.fitzp.sum}
g_1 \!\oplus\! g_2 (v,v^*) = g_1(v,v^*- \alpha_2(v)) + \langle \alpha_2(v),v\rangle
\qquad\forall (v,v^*)\in V \!\times\! V'.
\end{equation} 
\end{proposition}

\medskip
\noindent{\bf Proof.\/}  
(i) Let us set 
\begin{equation}
\varphi(v,v^*,z^*) = g_1(v,v^*- z^*) + g_2(v,z^*)
\qquad\forall (v,v^*,z^*)\in V \!\times\! V'\!\times\! V'.
\end{equation}
This auxiliary function is lower semicontinuous, and fulfills \eqref{eq.fitzp.coee} by \eqref{eq.fitzp.coo}.
By applying Lemma~\ref{lem1} with $X= V \!\times\! V'$, $Y= V'$, 
we then infer that $g_1 \!\oplus\! g_2$ is lower semicontinuous.

Let us check that $\varphi$ is convex. By the convexity of $g_1$ and $g_2$,  
for any $(v_i,v^*_i,z^*_i)\in V \!\times\! V'\!\times\! V'$ ($i=1,2$)
and any $\lambda\in [0,1]$,
\[
\varphi(\lambda(v_1,v^*_1,z^*_1) + (1-\lambda)(v_2,v^*_2,z^*_2))
\le\lambda\varphi(v_1,v^*_1,z^*_1) + (1-\lambda) \varphi(v_2,v^*_2,z^*_2). 
\]
By taking the infimum with respect to $z^*_1$ and $z^*_2$, we then get 
\[
g_1 \!\oplus\! g_2(\lambda(v_1,v^*_1)+ (1-\lambda)(v_2,v^*_2)) 
\le \lambda g_1 \!\oplus\! g_2(v_1,v^*_1)+ (1-\lambda)g_1 \!\oplus\! g_2(v_2,v^*_2);
\]
thus $g_1 \!\oplus\! g_2$ is convex. 
It is straightforward to check that $g_1 \!\oplus\! g_2\ge \pi$,
and that $g_1 \!\oplus\! g_2$ represents the operator $\alpha_1 + \alpha_2$.

(ii) If $\alpha_2:V\to V'$ is linear and continuous, then $\alpha_2$ is represented by 
$g_2: (v,v^*)\mapsto \langle v^*,v\rangle + I_{\alpha_2}(v,v^*)$.
\eqref{eq.fitzp.sum} then follows from \eqref{eq.fitzp.infcon}.
\hfill$\Box$

\bigskip 
\noindent{\bf Functional set-up.} 
First we introduce some spaces of time-dependent functions.
Let us assume that we are given a triplet of real Hilbert spaces
\footnote{ Here we deal with the Hilbert set-up and $p=2$ for technical reasons, 
that we shall see in the next section.
} 
\begin{equation}\label{eq.BEN.spaces1}
V\subset H=H'\subset V' \hbox{ \ with continuous and dense injections, }
\end{equation} 
whence
\begin{equation} 
\begin{split}
&L^2(0,T;V) \subset L^2(0,T;H) = L^2(0,T;H)' 
\subset L^2(0,T;V)' = L^2(0,T;V')
\\
&\hbox{with continuous and dense injections. } 
\end{split} 
\end{equation}
Let us also set
\begin{equation}\label{eq.BEN.spaces2}
{\cal V} := \big\{v\in L^2(0,T;V): 
D_tv \in L^2(0,T;V')\big\} \; (\subset C^0([0,T];H)),
\end{equation}
which, equipped with the graph norm, is a Hilbert space. 
We shall identify the spaces
\begin{equation} 
L^2(0,T;V)\!\times\! L^2(0,T;V') \simeq L^2(0,T; V\!\times\! V').
\end{equation}
 
\noindent{\bf An initial-value problem.}
Let us assume that
\begin{equation}\label{eq.BEN.mmon}
\alpha: V\to {\cal P}(V') \hbox{ is maximal monotone,} 
\end{equation}  
fix any $u^*\in L^2(0,T;V')$ and any $u^0\in H$, and consider the initial-value problem  
\begin{equation}\label{CauPb1}  
\left\{\begin{split} 
&u\in {\cal V},
\\
&D_tu + \alpha(u) \ni u^* \qquad \hbox{ in $V'$, a.e.\ in }{}]0,T[,
\\
&u(0) = u^0.
\end{split}\right.
\end{equation} 

We shall assume that
\begin{equation}\label{eq.BEN.ine}
\exists A,B>0: \forall v\in V, \forall v^*\in \alpha(v) 
\qquad
\|v^*\|_{V'} \le A\|v\|_V +B,
\end{equation} 
and define the global-in-time operator 
$\widehat\alpha: L^2(0,T;V)\to {\cal P}(L^2(0,T;V'))$: 
\begin{equation}\label{eq.BEN.param}
\begin{split}
&\forall (u,w)\in L^2(0,T;V)\times L^2(0,T;V'), \quad
\\
&w\in \widehat\alpha(u) 
\quad\Leftrightarrow\quad
w(t)\in \alpha(u(t)) \;\;\hbox{ for a.e.\ }t\in{} ]0,T[.
\end{split}
\end{equation}  
\indent 
The Cauchy problem~\eqref{CauPb1} is then equivalent to the following global-in-time problem
\begin{equation}\label{CauPb2}  
\left\{\begin{split} 
&u\in {\cal V},
\\
&D_tu + \widehat\alpha(u) \ni u^* \qquad\hbox{in }L^2(0,T;V'),
\\
&u(0) = u^0.
\end{split}\right.
\end{equation}

In the remainder of this section we shall deal with (proper) representable operators. 
As we saw, these are necessarily monotone, but need not be maximal monotone.

\begin{proposition}\label{timeint}
Let the assumptions \eqref{eq.BEN.mmon} and \eqref{eq.BEN.ine} be fulfilled.
Then: 
\\
\indent
(i) $\alpha$ is representable in ${\cal F}(V)$ if and only if $\widehat\alpha$
is representable in ${\cal F}(L^2(0,T;V))$.
More precisely, $\alpha$ is represented by a function $\varphi\in {\cal F}(V)$ 
if and only if $\widehat\alpha$ is represented by
$\Phi = \int_0^T \varphi(\cdot,\cdot) \, dt \in {\cal F}(L^2(0,T;V))$. 
\\
\indent
(ii) $\alpha$ is maximal monotone if and only if so is $\widehat\alpha$. 
\end{proposition}

\medskip
\noindent{\bf Proof.\/}   
(i) If $\varphi\in {\cal F}(V)$ then it is convex, therefore $\Phi$ is also convex and 
lower semicontinuous. 
Conversely, if $\Phi\in {\cal F}(V)$ then the same holds for $\varphi$, as it is easily checked by considering constant functions.

For any $(v,v^*)\in {\cal V}\!\times\! {\cal V}'$, 
\begin{equation}  
\begin{split}  
&\int_0^T \varphi(v(t),v^*(t)) \, dt \ge \int_0^T \langle v^*,v \rangle \, dt,
\\
&\varphi(v(t),v^*(t)) \ge \langle v^*(t),v(t)\rangle
\qquad\hbox{ for a.e.\ }t\in{} ]0,T[.
\end{split}  
\end{equation}  
If the first inequality is reduced to an equality, 
the same holds for the second inequality for a.e.\ $t$.
The converse implication is straightforward. Parts (i) is thus established.

(ii) Let us next represent the operator $\widehat{\alpha}$ by the function 
\[
\Phi = \int_0^T \varphi(\cdot,\cdot) \, dt: L^2(0,T;V \!\times\! V')
\to \erre\cup \{+\infty\}.
\] 
If $\alpha$ is maximal monotone, then by Theorem~\ref{teo.maxm} 
$\varphi^*\in {\cal F}(V')$ and it represents ${\alpha}^{-1}$. 
By part (i), $\Phi^*= \int_0^T \varphi^*(\cdot,\cdot) \, dt$ then represents 
the operator $\widehat{{\alpha}^{-1}}$.
As $\widehat{\alpha}^{-1} = \widehat{{\alpha}^{-1}}$,
again by Theorem~\ref{teo.maxm} we then infer that $\widehat{\alpha}$ is maximal monotone. The converse implication is straightforward.
\hfill$\Box$

\bigskip
\noindent{\bf Time-dependent operators.\/} 
If $\alpha $ explicitly depends on time, 
a precise formulation of the corresponding initial-value problem may be given by using
some classical notions of measurability, that we briefly recall.
A thorough account may be found e.g.\ in the monograph \cite{CaVa}. 
 
Let us denote by ${\cal B}(V)$ (${\cal L}(0,T)$, resp.) the 
$\sigma$-algebra of the Borel- (Lebesgue- resp.) measurable subsets of the 
separable space $V$ (the interval $]0,T[$, resp.), and by 
${\cal B}(V)\!\otimes\! {\cal L}(0,T)$ the product $\sigma$-algebra.
Let us assume that 
\begin{equation}\label{eq.BEN.tdep}
\begin{split}
&\alpha: V \!\times\! {}]0,T[{} \to {\cal P}(V') 
\hbox{ is a measurable and closed-valued proper operator.}
\end{split} 
\end{equation}
By this we mean that 
\[
\alpha^{-1}(A) := \big\{(v,t)\in V\!\times {}]0,T[{}: 
\alpha(v,t)\cap A\not= \emptyset \big\} \in {\cal B}(V)\!\otimes\! {\cal L}(0,T)
\]
for any open subset $A$ of $V'$, and that $\alpha(v,t)$ is closed for any $v\in V$ and a.e.\ $t$. 
For instance, this is fulfilled whenever $\alpha$ is maximal monotone. 
It is easy to see that, assuming \eqref{eq.BEN.tdep},
for any ${\cal L}(0,T)$-measurable mapping $v: {}]0,T[{} \to V$ the mapping 
${}]0,T[{} \to {\cal P}(V'): t\mapsto \alpha(v(t),t)$ is also measurable and closed-valued,
see e.g.\ Lemma~3.1 of \cite{Vi13}.
By a classical theorem, see e.g.\ Sect.~III.2 of \cite{CaVa}, this composed function has then a measurable selection; that is, there exists a measurable function
$w: {}]0,T[{}\to V'$ such that $w(t)\in \alpha(v(t),t)$ for a.e.\ $t$. 

\medskip
Under the assumption \eqref{eq.BEN.tdep},
we can now extend the initial-value problem~\eqref{CauPb1} to time-dependent operators:
\begin{equation}\label{CauPb1+}  
\left\{\begin{split} 
&u\in {\cal V},
\\
&D_tu + \alpha(u,t) \ni u^* \quad \hbox{ in $V'$, a.e.\ in }{}]0,T[,
\\
&u(0) = u^0.
\end{split}\right.
\end{equation}
This also is equivalent to a global-in-time formulation
analogous to problem \eqref{CauPb2}.

\begin{remarks}\rm
(i) The existence of measurable representative functions of time-dependent maximal monotone operators was addressed in Section~3 of \cite{Vi13}.

(ii) Rather then dealing with measurable multi-valued operators, 
it seems simpler to proceed via representing functions, as next we do. 
\end{remarks}
 
\noindent{\bf Extended BEN principle.}
Let us assume that \eqref{eq.BEN.tdep} is fulfilled and that
\begin{equation}\label{eq.BEN.repmes}
\begin{split}
&\hbox{for a.e.\ $t\in {} ]0,T[$, $\alpha(\cdot,t)$ is represented by 
$\varphi(\cdot,\cdot,t)\in {\cal F}(V)$, with }
\\
&\varphi: V\times V'\times {}]0,T[{}\to \erre\cup \{+\infty\} \ \ \
{\cal B}(V\times V')\!\otimes\! {\cal L}(0,T) \hbox{-measurable.}
\end{split} 
\end{equation} 
Let us fix any $u^0\in H$, define the affine subspace
\begin{equation}\label{eq.BEN.spaces3}
{\cal V}_{u^0} := \big\{v\in L^2(0,T;V): D_tv \in L^2(0,T;V'): v(0) = u^0\big\}
\subset {\cal V},
\end{equation}
and the nonnegative functional
\begin{equation}
\begin{split} 
\Phi(v,v^*) 
&:= \int_0^T [\varphi(v,v^* -D_tv,t) - \langle v^*- D_tv,v\rangle] \, dt   
\\
&= \int_0^T [\varphi(v,v^* -D_tv,t) - \langle v^*,v\rangle] \, dt 
+ {1\over2} \|v(T)\|_H^2 - {1\over2} \|u^0\|_H^2  
\\
& \qquad \qquad \qquad \qquad\forall (v,v^*)\in {\cal V}_{u^0} \!\times\! L^2(0,T;V'), 
\\
\Phi(v,v^*) 
&:= +\infty \qquad \qquad \hbox{ for any other }(v,v^*)\in {\cal V} \!\times\! {\cal V}'.
\end{split}
\label{eq.BEN.BENfun} 
\end{equation}

This definition is legitimate, since the function 
${\cal V}\to\erre: v\mapsto \|v(T)\|_H^2$ 
is well defined (and weakly lower semicontinuous),
as ${\cal V} \subset C^0([0,T];H)$ with continuous injection.
    
\begin{theorem}[Extended BEN principle] \label{teo.repr}
Assume that \eqref{eq.BEN.tdep} and \eqref{eq.BEN.repmes} are fulfilled, and define
$\widehat\alpha$ and $\Phi$ as in \eqref{eq.BEN.param} and \eqref{eq.BEN.BENfun},
for some $u^0\in H$. 
Then:
\\
\indent
(i) $\Phi\in {\cal F}({\cal V})$ and it represents the operator 
$D_t + \widehat\alpha: {\cal V}\to {\cal V}'$
associated with the initial condition $u(0) = u^0$. 
\\
\indent
(ii) For any $u^*\in V'$, the initial-value problem \eqref{CauPb1} is equivalent to the following null-minimization problem:
\begin{equation}\label{eq.BEN.min+} 
\left\{\begin{split} 
&u\in {\cal V},
\\[1mm]
&\Phi(u,u^*) = 0 \;\; \big(\! =\inf_{\cal V} \Phi(\cdot,u^*) \big),
\\
&u(0) = u^0.
\end{split}\right.
\end{equation}
\end{theorem}

This equivalence if what we shall properly refer to as the {\it extended BEN principle.\/}

\medskip
\noindent{\bf Proof.\/}  
Part (i) stems from part (ii) of Proposition~\ref{sum}, with $\alpha_1 = \widehat\alpha$,
$\alpha_2 = D_t$ and the space ${\cal V}$ in place of $V$.  
 
As $\varphi$ represents $\alpha$, for any $(v,v^*)\in {\cal V} \times L^2(0,T;V')$ 
the first integrand of \eqref{eq.BEN.BENfun} is nonnegative; 
moreover, it vanishes a.e.\ in $]0,T[$ if and only if $v$ solves \eqref{CauPb1}$_2$
and \eqref{CauPb1}$_3$.
This entails that $\Phi(v,v^*) \ge 0$ for any $(v,v^*)\in {\cal V} \times L^2(0,T;V')$,
and that the equation $\Phi(u,u^*) = 0$ is equivalent to \eqref{CauPb1}$_2$
and \eqref{CauPb1}$_3$.  
\hfill$\Box$

\begin{remark}\rm
Two issues may be addressed:

(i) Existence of a solution of the initial-value problem \eqref{CauPb1}
via the extended BEN principle Theorem~\ref{teo.repr}.

(ii) Structural properties of compactness and stability of the initial-value problem,
in the sense that we illustrate in the next section.
In this case we also let $u^*$ vary, and minimize $\Phi$ w.r.t.\ the pair $(u,u^*)$,
thus with no need of prescribing the minimum to vanish. 
\end{remark}

\section{$\Gamma$-compactness and $\Gamma$-stability of minimization problems} 
\label{sec.comstab}

\noindent
In this section we illustrate the structural compactness and stability of minimization principles. 
This will rest on $\Gamma$-convergence w.r.t.\ what we shall refer to as a 
{\it nonlinear weak topology.\/}

\bigskip
\noindent{\bf Structural compactness and stability of minimization problems.}  
Let $X$ be a topological space
and ${\cal G}$ be a family of functionals $X\to \erre\cup\{+\infty\}$,
equipped with a suitable topology, or at least a notion of convergence. 
We shall say that the problem of minimizing these functionals is 
{\it structurally compact\/} if the family ${\cal G}$ is sequentially compact, and
the corresponding minimizers are confined to a sequentially relatively compact subset of $X$.
\footnote{ We restrict ourselves to sequential concepts since this suffices for our purposes.
} 

We shall say that the minimization problem is {\it structurally stable\/} if 
\begin{equation}
\left\{\begin{split} 
&\Psi_n(u_n) - \inf \Psi_n \to 0
\\ 
&u_n\to u \quad\hbox{ in }X
\\ 
&\Psi_n \to \Psi\quad\hbox{ in }{\cal G}
\end{split} 
\qquad\Rightarrow\quad
\Psi(u) = \inf \Psi. \right.
\end{equation} 
(The condition $\Psi_n(u_n) - \inf \Psi_n \to 0$ is obviously more general than 
$\Psi_n(u_n) = \inf \Psi_n$ for any $n$.)
These definitions also apply to null-minimization problems.

The selection of the topology of the family of functionals ${\cal G}$ is crucial, 
and may not be an obvious choice. 
Structural compactness and stability are in competition:
the topology must be sufficiently weak in order to allow for compactness,
and at the same time it must be so strong to provide stability
(essentially, passage to the limit in the perturbed minimization problem). 
The exigence of compactness suggests the use of a weak-type topology.
We shall see that here $\Gamma$-convergence is especially appropriate, 
more than other notions of variational convergence like the Mosco-convergence
(see e.g.\ \cite{At}, \cite{Mo}), namely the simultaneous $\Gamma w$- and $\Gamma s$-convergence to a same function.

\bigskip
\noindent{\bf On weak $\Gamma$-convergence.\/}
Let the spaces $H, V, {\cal V}$ be defined as in \eqref{eq.BEN.spaces1} and 
\eqref{eq.BEN.spaces2}. 
For any $n$, let $\varphi_n\in {\cal F}(V)$ represent a proper measurable operator 
$\alpha_n: V\to {\cal P}(V')$, and define the functionals $J_n$ and $\Phi_n$ as in
\eqref{eq.fitzp.nullmindef} and \eqref{eq.BEN.BENfun}:
\begin{eqnarray} 
&J_n(v,v^*) := \varphi_n(v,v^*) - \langle v^*,v\rangle 
\qquad\forall (v,v^*)\in V \!\times\! V',
\label{eq.comstab.motiv1}
\\
&\begin{split}
\Phi_n(v,v^*) := \int_0^T \varphi_n(v,v^* -D_tv) dt 
- \int_0^T \langle v^*-D_tv,v\rangle \, dt &
\\
\forall (v,v^*)\in {\cal V} \!\times\! L^2(0,T;V'). &
\end{split} 
\label{eq.comstab.motiv2} 
\end{eqnarray} 

Let us first consider the (sequential) $\Gamma$-convergence of the sequence $\{J_n\}$
in the weak topology of $V \!\times\! V'$. 
This requires passing to the inferior limit in $\langle v_n^*,v_n\rangle$, as
$v_n\wto v$ in $V$ and $v_n^*\wto v^*$ in $V'$.
This difficulty may be removed simply by confining $\{v_n^*\}$, e.g., to a bounded subset of 
$H$, if
\begin{equation}\label{eq.comstab.triple}
V\subset H=H' \subset V' 
\hbox{ \ \ with continuous, dense and compact injections.}
\end{equation} 
 
On the other hand, dealing with the $\Gamma$-convergence of the sequence 
$\{\Phi_n\}$, one must pass to the limit in the integral 
\begin{equation}\label{eq.comstab.int} 
- \int_0^T \! \langle v_n^*,v_n\rangle \, dt + 
\int_0^T \! \langle D_tv_n,v_n\rangle \, dt
= \int_0^T \! \langle v_n^*,v_n\rangle \, dt 
+ {1\over2} \|v_n(T)\|_H^2 - {1\over2} \|u_n(0)\|_H^2,
\end{equation} 
and this looks less obvious than for $\{J_n\}$.

\medskip 
\noindent{\bf A different functional set-up.} 
Let us define the measure
\begin{equation}\label{eq.comstab.mu'}
\mu(A) = \int_A (T-t) \, dt \quad\forall A\in {\cal L}(0,T),
\quad\hbox{ i.e., }\quad
d\mu(t) = (T-t) dt,
\end{equation}
introduce the weighted spaces
\begin{equation}\label{eq.comstab.spaces4} 
\begin{split} 
&L^2_\mu(0,T;H) := \Big\{\mu\hbox{-measurable } v:{} ]0,T[{}\to H :
\!\! \int_0^T \!\|v(t)\|_H^2 \, d\mu(t) <+\infty \Big\}, 
\\
&L^2_\mu(0,T;V) := \Big\{v\in L^2_\mu(0,T;H): 
\!\int_0^T \|v(t)\|_V^2 \, d\mu(t) <+\infty \Big\},
\\
&{\cal V}_\mu := \big\{v\in L^2_\mu(0,T;V): D_tv \in L^2_\mu(0,T;V')\big\},
\\ 
\end{split} 
\end{equation}
and equip them with the respective graph norm.
\footnote{ The reason for introducing the weight $T-t$ will appear ahead.
}
These are Hilbert spaces, and by identifying $L^2_\mu(0,T;H)$ with its dual space
we get  
\begin{equation} 
\begin{split} 
&L^2_\mu(0,T;V') = \Big\{\mu\hbox{-measurable } v^*:{} ]0,T[{}\to V' : 
\!\! \int_0^T \! \|v^*(t)\|_{V'}^2 \, d\mu(t) <+\infty \Big\},
\\
&L^2_\mu(0,T;V) \subset L^2_\mu(0,T;H) = L^2_\mu(0,T;H)' 
\subset L^2_\mu(0,T;V)' = L^2_\mu(0,T;V')
\\
&\hbox{with continuous and dense injections. } 
\end{split} 
\end{equation}
 
Notice that ${\cal V}_\mu \subset C^0([0,T[;H)$.
Under the assumption \eqref{eq.comstab.triple}, 
\begin{equation}\label{eq.comstab.y0'}
\hbox{ the injection ${\cal V}\to L^2(0,T;H)$ 
$(\subset L^2_\mu(0,T;H))$ is compact, } 
\end{equation} 
at variance with the injections 
${\cal V}\to C^0([0,T];H)$ and ${\cal V}_\mu\to L^2(0,T;H)$.
If $v_n\wto v$ in ${\cal V}$ thus $v_n\to v$ in $L^2(0,T;H)$, so that
one can pass to the limit in the first integral of \eqref{eq.comstab.int}, provided that
$v_n^*\wto v^*$ in $L^2(0,T;H)$. 
But we are left with the problem of passing to the limit in the second integral of the same formula, and this will be achieved via a modification of the functional.

For any $n$, let us now replace $dt$ by $d\mu(t)$ in the definition of the functional 
$\Phi_n$, see \eqref{eq.comstab.motiv1}, and set
\begin{eqnarray} 
&\begin{split}
\Phi_n(v,v^*) := \int_0^T \varphi_n(v,v^* -D_tv) \, d\mu(t) 
- \int_0^T \langle v^*-D_tv,v\rangle \, d\mu(t) &
\\
\forall (v,v^*)\in {\cal V} \!\times\! L^2(0,T;V'). &
\end{split} 
\label{eq.comstab.motiv2=} 
\end{eqnarray} 
(As we shall see, this corresponds to applying a double time integration to the time-integrand
of \eqref{eq.comstab.motiv2}.)

We must then pass to the limit in the integral 
\begin{equation}\label{eq.comstab.dueint}  
- \int_0^T \! \langle v_n^*,v_n\rangle \, d\mu(t) + 
\int_0^T \! \langle D_tv_n,v_n\rangle \, d\mu(t).
\end{equation}
Concerning the second integral, assuming that $v_n(0)\to v(0)$ in $H$, we have
\begin{equation}\label{eq.comstab.motiv3}
\begin{split}
&\int_0^T \! \langle D_tv_n,v_n\rangle \, d\mu(t) 
= {1\over2} \int_0^T D_t \big(\|v_n(t)\|_H^2\big) \, (T-t) \, dt 
\\
&= {1\over2} \int_0^T \! \|v_n\|_H^2 \, dt - {T\over2} \|v_n(0)\|_H^2 
\\
&\to {1\over2} \int_0^T \! \|v\|_H^2 \, dt  - {T\over2} \|v(0)\|_H^2 
= \!\int_0^T \! \langle D_tv,v\rangle \, d\mu(t).
\end{split}
\end{equation} 

\begin{remarks}\rm
(i) The integral functional
$v\mapsto \int_0^T \! \langle D_tv,v\rangle \, d\mu(t)$ is weakly continuous on 
${\cal V}$. The functional $v\mapsto \int_0^T \! \langle D_tv,v\rangle \, dt$
instead is just weakly lower semicontinuous on the same space; 
this entails the existence of a minimizer, but not the $\Gamma$-convergence.
This is the main reason why we introduced the weight $T-t$.

(ii) On the other hand, if one prescribes the solution $v$ to be $T$-periodic in time,
then $\int_0^T \langle D_tv_n,v_n\rangle \, dt =0$ for any $n$.
In this case it is not needed to introduce any weight function.
\end{remarks}
 
\noindent{\bf A nonlinear topology of weak type.\/}
The convergence of the terms
\[
\langle v_n^*,v_n\rangle
\qquad\hbox{ and }\qquad
\int_0^T \langle v_n^*- D_tv_n,v_n\rangle \, d\mu(t)
\] 
yields the inferior-limit condition of the definition of $\Gamma$-convergence
of the sequences $\{J_n\}$ and $\{\widetilde\Phi_n\}$.
But one must also pass to the limit on a so-called recovery sequence.
This prompts us to complement the weak topology of $V\!\times\! V'$ 
with the convergence $\langle v_n^*,v_n \rangle \to \langle v^*,v\rangle$
when dealing with the sequence $\{J_n\}$,
and similarly to complement the weak topology of 
\begin{equation}
L^2_\mu(0,T;V \!\times\! V') \simeq
L^2_\mu(0,T;V) \!\times\! L^2_\mu(0,T;V')
\end{equation}
with the convergence 
$\int_0^T \langle v_n^*,v_n \rangle \, d\mu(t) 
\to \int_0^T \langle v^*,v\rangle \, d\mu(t)$ 
when dealing with $\{\widetilde\Phi_n\}$.

More specifically, we name {\it nonlinear weak topology\/} of $V\!\times\! V'$, 
and denote by $\widetilde\pi$, the coarsest among the topologies of this space  
that are finer than the product of the weak topology of $V$ by the weak topology of $V'$, 
and for which the mapping $\pi$ (that we defined in \eqref{eq.fitzp.pi}) is continuous.
For any sequence $\{(v_n,v^*_n)\}$ in $V\!\times\! V'$, thus 
\begin{equation}
\begin{split}
&(v_n,v_n^*)\pito (v,v^*) \quad\hbox{in }V\!\times\! V'Ê
\quad\Leftrightarrow 
\\
&v_n\wto v \hbox{ \ \ in }V, \quad
v_n^*\wto v^* \hbox{ \ \  in }V', \quad
\langle v_n^*,v_n \rangle \to \langle v^*,v\rangle,
\end{split}
\end{equation}
and similarly for nets. (The nonlinearity is obvious: 
a linear combination of two converging sequences need not converge.)

This construction is extended to the space 
$L^2_\mu(0,T;V \!\times\! V')$ 
in an obvious way: in this case the duality mapping reads 
$L^2_\mu(0,T;V \!\times\! V')\to \erre: 
(v,v^*)\mapsto \int_0^T \langle v^*,v\rangle \, d\mu(t)$, and we set
\begin{equation}
\begin{split}
&(v_n,v_n^*)\pito (v,v^*) 
\quad\hbox{in }L^2_\mu(0,T;V \!\times\! V')Ê
\quad\Leftrightarrow 
\\
&v_n\wto v \hbox{ \ in }L^2_\mu(0,T;V), \;
v_n^*\wto v^* \hbox{ \  in }L^2_\mu(0,T;V') \; \hbox{ and }
\\
&\int_0^T \! \langle v_n^*,v_n \rangle \, d\mu(t) \to 
\! \int_0^T \! \langle v^*,v\rangle \, d\mu(t)
\end{split}
\end{equation}
and similarly for nets.

\bigskip
\noindent{\bf {\boldmath $\Gamma\widetilde\pi$}-compactness and 
{\boldmath $\Gamma\widetilde\pi$}-stability of ${\cal F}(V)$.\/}
As the weak topology and the nonlinear weak-type topology $\widetilde\pi$ are nonmetrizable, some caution is needed in dealing with {\it sequential\/} 
$\Gamma$-convergence w.r.t.\ either topology,
see e.g.\ \cite{At}, \cite{Da}.
For functions defined on a topological space, the definition of 
$\Gamma$-convergence involves the filter of the neighborhoods of each point.
If the space is metrizable, that notion may equivalently be formulated in terms of the family of converging sequences, but this does not hold in general. 
We shall refer to these two notions as {\it topological\/} and {\it sequential\/} $\Gamma$-convergence, respectively.
If not otherwise specified, reference to the topological notion should be understood.
 
It is known that bounded subsets of a separable and reflexive space equipped with the weak topology are metrizable.
The same holds for the nonlinear weak topology $\widetilde\pi$ of $V\!\times\! V'$, 
see Section~4 of \cite{ViCalVar}. 
This property is at the basis of the next statement, which gathers $\Gamma$-compactness, $\Gamma$-closedness of ${\cal F}(V)$, and $\Gamma$-stability. 
We shall denote by $\gr(\alpha)$ the graph of an operator $\alpha$, i.e.,
$\gr(\alpha) = \{(v,v^*): v^*\in \alpha(v)\}$.

\begin{theorem} [\cite{ViCalVar}] \label{comsta} 
Let $V$ be a separable real Banach space, and 
$\{\psi_n\}$ be an equi-coercive sequence in ${\cal F}(V)$, in the sense that 
\begin{equation}\label{eq.equicoer}
\sup_{n\in\enne} \big\{\|v\|_V + \|v^*\|_{V'}: 
(v,v^*)\in V\!\times\! V', \psi_n(v,v^*) \le C \big\}<+\infty
\quad\forall C\in \erre.
\end{equation}
\indent 
Then:  
(i) There exists $\psi: V\!\times\! V' \to \erre\cup\{+\infty\}$ such that, possibly extracting a subsequence, $\psi_n$ $\Gamma\widetilde\pi$-converges to $\psi$ both topologically and sequentially.
\hfill\break\indent
(ii) This entails that $\psi \in {\cal F}(V)$.
\hfill\break\indent
(iii) If $\alpha_n$ ($\alpha$, resp.) is the operator that is represented by 
$\psi_n$ ($\psi$, resp.) for any $n$, then 
$K\tau \limsup_{n\to \infty} \gr(\alpha_n) \subset \gr(\alpha)$
in the sense of Kuratowski, i.e., 
\begin{equation}\label{eq.comstab.lim}
\forall \text{ sequence }\{(v_n,v_n^*)\in \gr(\alpha_n)\}, \quad
(v_n,v^*_n)\pito (v,v^*)
\quad\Rightarrow\quad v^*\in \alpha(v). 
\end{equation}
\indent
(iv) The same results hold if $V$ is replaced by $L^2_\mu(0,T;V)$ throughout.
\end{theorem}  

\medskip
\noindent{\bf Proof.\/}  
Part (i) is Theorem~4.4 of \cite{ViCalVar}.
 
The $\widetilde\pi$-lower semicontinuity of the $\psi_n$s and the property 
``$\psi_n \ge \pi$'' (see \eqref{eq.fitzp.pi}) are preserved by  
$\Gamma\widetilde\pi$-convergence.
If $\psi_n \in {\cal F}(V)$ for any $n$, namely the $\psi_n$s are also convex,
then the same holds for $\psi$, which is then lower semicontinuous.
Part (ii) is thus established. 

Next let $\alpha_n$ ($\alpha$, resp.) and the sequence $\{(v_n,v^*_n)\}$ be as prescribed 
in part (iii). 
Recalling the inferior-limit condition of the definition of $\Gamma\widetilde\pi$-convergence
and the definition of $\widetilde\pi$-convergence, we have
\begin{equation}\label{eq.comstab.lim'}
\psi(v,v^*) 
\le\liminf_{n\to \infty} \psi_n(v_n, v^*_n)
\overset{\eqref{eq.fitzp.convrep}_2}{\le}
\liminf_{n\to \infty} \; \langle v^*_n,v_n \rangle = \langle v^*,v \rangle.
\end{equation}
Thus $v^*\in \alpha(v)$, as $\psi$ represents $\alpha$. 
The implication \eqref{eq.comstab.lim} is thus established.  

Part (iv) is easily checked by mimicking the above procedure.
\hfill$\Box$
 
\begin{remark}\rm 
In general $\alpha_n$ does not converge in the sense of Kuratowski.
For instance, defining $f_b$ as in \eqref{eq.fitzp.counterex} for any $b>0$, 
if $b_n\to 1/2$ then $f_{b_n}$ $\Gamma\widetilde\pi$-converges to $f_{1/2}$, 
but the represented operators $\alpha_n$ do not converge in the sense of Kuratowski
--- actually, $\gr(\alpha)\not\subset K\tau \liminf_{n\to \infty} \gr(\alpha_n)$.  
\end{remark}  
 
After finding out that the class ${\cal F}(L^2_\mu(0,T;V))$ is stable by 
$\Gamma\widetilde\pi$-convergence, the next question arises.
Let $\{\varphi_n\}$ be a sequence of representative functions of ${\cal F}(V)$, and define the memoryless superposition operators
\begin{equation}
\psi_n: L^2_\mu(0,T;V \!\times\! V')\to L^1(0,T): w\mapsto \varphi(w) \quad\forall n.
\end{equation} 
If $\psi_n$ $\Gamma\widetilde\pi$-converges to some operator $\psi$ in the sense of 
\eqref{eq.evol.tesi}, is then $\psi$ necessarily memoryless, too? 
A positive answer is provided by the next statement, 
which is essentially a particular case of Theorem~\ref{teo.comp}, 
augmented with the hypothesis \eqref{eq.comstab.reppn} and with the corresponding property 
in the limit.
\footnote{ Here we deal with the Hilbert set-up and $p=2$ for technical reasons, 
that we shall see in the next section.
} 
Here we also allow representing functions to depend on time. 

\begin{theorem} \label{teo.comp'}
Let $V$ be a real separable Hilbert space,
and $\mu$ be the measure on $]0,T[$ such that $d\mu(t) = (T-t) \, dt$.
Let $\{\varphi_n\}$ be a sequence of normal functions 
${}]0,T[{} \times V \!\times\! V'\to \erre^+$ such that  
\begin{eqnarray} 
&\varphi_n(t,\cdot) \in {\cal F}(V)
\qquad\hbox{for a.e.\ }t\in {}]0,T[, \forall n,  
\label{eq.comstab.reppn}
\\
&\begin{split}
&\exists C_1,C_2,C_3 >0: \forall n,\hbox{for a.e.\ }t\in {}]0,T[, \forall w\in V \!\times\! V', 
\\
&C_1 \|w\|_{V \!\times\! V'}^2\le \varphi_n(t,w) \le C_2\|w\|_{V \!\times\! V'}^2 +C_3, 
\end{split}   
\label{eq.comstab.equibc'}  
\\
&\varphi_n(t,0) =0 \qquad\hbox{ for a.e.\ }t\in {}]0,T[,\forall n,
\label{eq.comstab.nul'}
\end{eqnarray} 
and define the operators $\psi_n: L^2_\mu(0,T;V \!\times\! V')\to L^1_\mu(0,T)$ by
\begin{equation}\label{eq.comstab.super'} 
\psi_{n,w}(t) = \varphi_n(t,w(t))  
\qquad\forall w\in L^2_\mu(0,T;V \!\times\! V'),\hbox{ for a.e.\ }t\in {}]0,T[,\forall n.
\end{equation}  
\indent
Then there exists a normal function $\varphi: {}]0,T[{} \times V \!\times\! V'\to \erre^+$ 
such that
\begin{equation}\label{eq.comstab.repp}
\varphi(t,\cdot) \in {\cal F}(V)
\qquad\hbox{for a.e.\ }t\in {}]0,T[,
\end{equation}  
and such that, defining the corresponding operator 
$\psi: L^2_\mu(0,T;V \!\times\! V')\to L^1_\mu(0,T)$ 
as in \eqref{eq.comstab.super'},
\begin{equation}\label{eq.comstab.tesi'} 
\begin{split}
&\hbox{ $\psi_n$ sequentially $\Gamma$-converges to }\psi
\\
&\hbox{ in the topology $\widetilde\pi$ of $L^2_\mu(0,T;V \!\times\! V')$ and } 
\\
&\hbox{ in the weak topology of $L^1_\mu(0,T)$ (cf.\ \eqref{eq.evol.defgamma.2}). }
\end{split}
\end{equation}
\indent
Moreover, if $\varphi_n$ does not explicitly depend on $t$ for any $n$, 
then the same holds for $\varphi$.
\end{theorem}

\medskip
\noindent{\bf Proof.\/}
Let us apply Theorem~\ref{teo.comp} with $X= V \!\times\! V'$,
$p=2$ and the topology $\tau = \widetilde\pi$; the hypothesis \eqref{eq.evol.gamcom} 
is indeed fulfilled for this topology, because of Theorem~4.4 of \cite{ViCalVar}. 
It then suffices to show that \eqref{eq.comstab.reppn} entails \eqref{eq.comstab.repp}.

The lower semicontinuity of $\varphi(t,\cdot)$ for a.e.\ $t$ is implicit in the normality of 
$\varphi$. In order to prove \eqref{eq.fitzp.convrep} for the function $\varphi$, let us set
\begin{equation}\label{eq.comstab.nonconv=}
\begin{split}
&J_n(t,v,v^*) := \varphi_n(t,v,v^*) - \langle v^*,v\rangle
\\[1mm]
&J(t,v,v^*) := \varphi(t,v,v^*) - \langle v^*,v\rangle 
\end{split}
\qquad\forall (v,v^*)\in V \!\times\! V', \hbox{ for a.e.\ }t,\forall n.
\end{equation}  
By \eqref{eq.comstab.reppn}, 
\begin{equation}\label{eq.comstab.nonconv+}
\int_A J_n(t,v,v^*) \, d\mu(t) \ge 0 
\qquad\forall (v,v^*)\in V \!\times\! V', \forall A\in {\cal L}(0,T),
\end{equation} 
and by \eqref{eq.comstab.tesi'} this inequality is preserved in the limit. 
Therefore $J(t,v,v^*)\ge 0$ a.e.\ in $]0,T[$.
As $\varphi$ is convex and lower semicontinuous, \eqref{eq.comstab.repp} follows.
\hfill$\Box$ 

\bigskip
\noindent{\bf Time-integrated (extended) BEN principle.}
Let us first notice that
\begin{equation}\label{eq.comstab.double} 
\begin{split}
&\int_0^T \langle D_tv,v\rangle \, d\mu(t) 
\overset{\eqref{eq.comstab.mu'}}{=} 
\int_0^T \! d\tau\! \int_0^\tau \! \langle D_tv,v\rangle \, dt 
\\
&= {1\over2}\! \int_0^T\! d\tau\! \int_0^\tau \! D_t\big(\|v(t)\|_H^2\big) \, dt 
= {1\over2} \!\int_0^T \! \|v(\tau)\|_H^2 \, d\tau - {T\over2} \|v(0)\|_H^2
\qquad\forall v\in {\cal V}.
\end{split}
\end{equation}
Let us fix any $u^0\in H$, and define the affine subspace
\begin{equation}\label{eq.comstab.spaces3}
{\cal V}_{\mu,u^0} := \big\{v\in L^2_\mu(0,T;V): 
D_tv \in L^2_\mu(0,T;V'): v(0) = u^0\big\}.
\end{equation}

Let $\varphi\in {\cal F}(V)$ represent a proper measurable operator 
$\alpha: V\to {\cal P}(V')$.
By a double time-integration, let us then define the nonnegative functional  
\begin{equation} \label{eq.comstab.BENfun'} 
\begin{split}
&\widetilde\Phi(v,v^*) 
:= \int_0^T \! d\tau\! \int_0^\tau \!\big[\varphi(v,v^* -D_tv) 
- \langle v^*-D_tv,v\rangle\big] \, dt  
\\
&\qquad\quad\overset{\eqref{eq.comstab.double}}{=}
\!\!\! \int_0^T \!\!\! \big[\varphi(v,v^* -D_tv)- \langle v^*,v\rangle\big] \, d\mu(t)
+{1\over2} \int_0^T \!\! \|v(t)\|_H^2 \, dt - {T\over2} \|u^0\|_H^2
\\
&\qquad\qquad\qquad\qquad\quad
\qquad\forall (v,v^*)\in {\cal V}_{\mu,u^0} \!\times\! L^2(0,T;V'), 
\\
&\widetilde\Phi(v,v^*) 
:= +\infty \qquad \qquad \hbox{ for any other }(v,v^*)\in {\cal V} \!\times\! {\cal V}'.
\end{split}
\end{equation}

The density of ${\cal V}\cap {\cal V}_\mu$ in ${\cal V}_\mu$ yields the next statement.

\begin{theorem}[Time-integrated (extended) BEN principle] \label{teo.repr'}
Let a proper operator $\alpha: V\to {\cal P}(V')$ be represented by a function 
$\varphi \in {\cal F}(V)$, and fulfill the condition \eqref{eq.BEN.ine}.
Let $\widehat\alpha$ and $\widetilde\Phi$ be defined as in \eqref{eq.BEN.param} and 
\eqref{eq.comstab.BENfun'}. 
\\
\indent
For any $u^*\in L^2_\mu(0,T;V')$ and $u^0\in H$, the initial-value problem \eqref{CauPb1} 
is then equivalent to the following null-minimization problem:
\begin{equation} 
u\in {\cal V},
\qquad
\widetilde\Phi(u,u^*) = 0 \;\; \big(\! =\inf_{\cal V} \widetilde\Phi(\cdot,u^*) \big).
\label{BEN'}
\end{equation}
\end{theorem}
  
\noindent{\bf Proof.\/} 
It suffices to notice that, as $\varphi\in {\cal F}(V)$, 
\eqref{CauPb1}$_2$ is equivalent to the minimization 
of the functional that is obtained by integrating the nonnegative function
$\varphi(v,v^* -D_tv) - \langle v^*-D_tv,v\rangle$ any number of times in $]0,T[$.
\hfill$\Box$

\section{Structural properties of maximal monotone flows}
\label{sec.par}

\noindent
In this section we show the structural compactness and stability of a class of 
flows governed by a maximal monotone operator.
We improve the results of \cite{ViCalVar}, on the basis of a different functional set-up 
and of Theorem~\ref{teo.comp}.
 
\bigskip
\noindent{\bf Quasilinear parabolic operators in abstract form.\/}  
Let $V$ and $H$ be real separable Hilbert spaces, and
\begin{equation}\label{eq.par.tripl}
V\subset H=H'\subset V' \hbox{ \ with continuous, dense and compact injections.} 
\end{equation} 
Let us define the spaces $L^2(0,T;H), L^2(0,T;V)$ and ${\cal V}$ as in \eqref{eq.BEN.spaces2}, 
and notice that in this case
\begin{equation}\label{eq.par.com} 
{\cal V} \subset L^2(0,T;H)
\hbox{ \ with continuous, dense and compact injection. } 
\end{equation}
Let us still define the measure $\mu$ as in \eqref{eq.comstab.mu'},
and identify the spaces
\begin{equation} 
L^2_\mu(0,T;V)\!\times\! L^2_\mu(0,T;V') \simeq L^2_\mu(0,T;V \!\times\! V').
\end{equation}

Let three sequences $\{\alpha_n\}, \{u^0_n\}$ and $\{h_n\}$ be given such that 
\begin{eqnarray}
&\forall n, \; \alpha_n: V\to {\cal P}(V') \hbox{ is maximal monotone, } 
\label{eq.par.alpha1}
\\
&\exists C_1, C_2>0: \forall n, 
\forall (v,v^*)\in \graph(\alpha_n), \quad 
\langle v^*,v \rangle \ge C_1|v\|_V^2 - C_2, 
\label{eq.par.alpha2}
\\ 
&\exists C_3, C_4>0: \forall n, 
\forall (v,v^*)\in \graph(\alpha_n), \quad 
\|v^*\|_{V'} \le C_3 \|v\|_V + C_4,  
\label{eq.par.alpha3}
\\
&\alpha_n(0) \ni 0 \qquad\forall n,
\label{eq.par.alpha3'}
\\
&u^0_n\to u^0 \quad\hbox{ in }H,
\\
&h_n\to h \quad\hbox{ in }L^2(0,T;V').
\label{eq.par.alpha4}
\end{eqnarray} 
The condition \eqref{eq.par.alpha3'} is not really restrictive: if it is not satisfied,
it may be recovered by selecting any $v^*_n\in \alpha_n(0)$ for any $n$, 
and then replacing $\alpha_n$ by $\bar\alpha_n(\cdot) = \alpha_n(\cdot) - v^*_n$. 
 
For any $n$, let us consider the following initial-value problem 
\begin{equation}\label{pbn1}  
\left\{\begin{split} 
&u_n\in {\cal V},
\\
&D_tu_n + \alpha(u_n) \ni h_n \qquad \hbox{ in $V'$, a.e.\ in }{}]0,T[,
\\
&u_n(0) = u^0_n.
\end{split}\right.
\end{equation} 
  
Next we recall a classical result of existence, uniqueness and boundedness.

\begin{lemma} [\cite{Bar}, \cite{Br1}, \cite{Ze1}] \label{regul}
Under the hypotheses above, for any $n$ the initial-value problem~\eqref{pbn1} has one and only one solution $u_n\in {\cal V}$. Moreover, the sequence $\{u_n\}$ is bounded in ${\cal V}$.
\end{lemma} 
 
\begin{remark}\rm 
Under strengthened assumptions on the data,
several regularity results are known to hold for the initial-value problem~\eqref{pbn1},
see e.g.\ \cite{Br1}. At variance with \cite{ViCalVar} here we do not use them,  
since we want to develop a method that rests on minimal assumptions,
in order to be able to extend it to more general flows, including:

(i) flows with a time-dependent (maximal monotone) operator,

(ii) doubly-nonlinear equations, that we deal with in the next two sections,

(iii) pseudo-monotone flows (to be addressed in a future work).
\end{remark}  

\begin{theorem}[Structural compactness and stability]\label{eq.par.stru}
Let \eqref{eq.par.tripl}--\eqref{eq.par.alpha4} be fulfilled, 
and for any $n$ let $u_n$ be the solution of problem~\eqref{pbn1}. Then:
\\
\indent
(i) There exists $u\in {\cal V}$ such that, possibly extracting a subsequence,  
\begin{equation}\label{eq.par.1a}
u_n\wto u \qquad\hbox{ in }{\cal V}.
\end{equation} 
\indent
(ii) There exists a function $\varphi\in {\cal F}(V)$ such that, setting 
\begin{eqnarray} 
&\varphi_n = (\pi+ I_{\alpha_n})^{**} \; (\in{\cal F}(V)),
\label{eq.par.de'}
\\[1mm]
&\begin{split} 
&\psi_{n,w}(t) = \varphi_n(w(t)), \qquad
\psi_w(t) = \varphi(w(t)),  
\\[1mm]
&\hbox{for a.e.\ }t\in{}]0,T[,\forall w\in L^2_\mu(0,T;V \!\times\! V'), \forall n,
\end{split} 
\label{eq.par.de''}
\end{eqnarray}
then $\psi_n,\psi: L^2_\mu(0,T;V \!\times\! V') \to L^1_\mu(0,T)$ and, 
possibly extracting a subsequence, 
\begin{equation}\label{eq.par.tesi+}
\begin{split}
&\hbox{ $\psi_n$ sequentially $\Gamma$-converges to }\psi
\\
&\hbox{ in the topology $\widetilde\pi$ of $L^2_\mu(0,T;V \!\times\! V')$ and } 
\\
&\hbox{ in the weak topology of $L^1_\mu(0,T)$ (cf.\ \eqref{eq.evol.defgamma.2}). }
\end{split}
\end{equation} 
\indent
(iii) Denoting by $\alpha: V\to {\cal P}(V')$ 
the monotone operator that is represented by $\varphi$,
$u$ solves the corresponding initial-value problem
\begin{equation}\label{pb=}  
\left\{\begin{split} 
&u\in {\cal V},
\\
&D_tu + \alpha(u) \ni h \qquad \hbox{ in $V'$, a.e.\ in }{}]0,T[,
\\
&u(0) = u^0.
\end{split}\right.
\end{equation}
\end{theorem}

\medskip
\noindent{\bf Proof.\/} 
For any $n$, let us define the affine subspace ${\cal V}_{\mu,u_n^0}$ of ${\cal V}$
as in \eqref{eq.comstab.spaces3}, and the nonnegative twice-time-integrated functionals  
as in \eqref{eq.comstab.BENfun'}:
\begin{equation} \label{eq.par.BENfun''} 
\begin{split}
&\widetilde\Phi_n(v,v^*) 
:= \int_0^T \! d\tau\! \int_0^\tau \big[\varphi_n(v,v^* -D_tv) - 
\langle v^*-D_tv,v\rangle\big] \, dt 
\\
&\overset{\eqref{eq.comstab.double}}{=}
\!\! \int_0^T \! \big[\varphi_n(v,v^* -D_tv) - \langle v^*,v\rangle\big] \, d\mu(t)
+{1\over2} \int_0^T \! \|v(t)\|_H^2 \, dt - {T\over2} \|u^0_n\|_H^2
\\
&\qquad\qquad\qquad \qquad\quad\; 
\qquad\forall (v,v^*)\in {\cal V}_{\mu,u^0_n} \!\times\! L^2_\mu(0,T;V'), 
\\
&\widetilde\Phi_n(v,v^*) 
:= +\infty \qquad \qquad \hbox{ for any other }(v,v^*)\in {\cal V} \!\times\! {\cal V}'.
\end{split}
\end{equation}

For any $n$, by Theorem~\ref{teo.repr'}, $u_n$ solves 
the initial-value problem~\eqref{pbn1} if and only if  
\begin{equation}\label{eq.par.min1}  
\left\{\begin{split} 
&u_n\in {\cal V}, 
\\[1mm]
&\widetilde\Phi_n(u_n,h_n) = 0\;\; 
\big(\! =\inf_{\cal V} \widetilde\Phi_n(\cdot,h_n) \big),
\\
&u_n(0) = u^0_n.
\end{split}\right.
\end{equation} 

By Lemma~\ref{regul}, \eqref{eq.par.1a} holds up to extracting a subsequence.
By \eqref{eq.comstab.motiv3} then
\begin{equation}\label{eq.par.co1}
\int_0^T \langle h_n - D_tu_n, u_n\rangle \, d\mu(t) \to
\int_0^T \langle h - D_tu, u\rangle \, d\mu(t),
\end{equation}
so that
\begin{equation}\label{eq.par.co2}
(h_n - D_tu_n, u_n)\pito (h - D_tu, u)
\qquad\hbox{ in }L^2_\mu(0,T;V \!\times\! V').
\end{equation}

Because of \eqref{eq.par.alpha2} and \eqref{eq.par.alpha3},
we may apply Theorem~\ref{teo.comp'}.
As here the functions $\varphi_n$'s are translation-invariant, 
there exists a function $\varphi\in {\cal F}(V)$ such that,
defining $\psi$ as in \eqref{eq.par.de''}, 
$\psi: L^2_\mu(0,T;V \!\times\! V') \to L^1_\mu(0,T)$ and, 
possibly extracting a subsequence, \eqref{eq.par.tesi+} is fulfilled. 
This yields
\begin{equation}\label{eq.par.ga} 
\begin{split}
F_n(v,v^*) := \int_0^T \! \varphi_n(v,v^*) \, d\mu(t) 
\Gammatopi \int_0^T \! \varphi(v,v^*) \, d\mu(t) =: F(v,v^*)
\\
\qquad\hbox{ sequentially in }V \!\times\! V'.
\end{split}
\end{equation}
 
Therefore
\begin{equation}
\begin{split}
&F(u, h- D_tu) 
\overset{\eqref{eq.par.co2},\eqref{eq.par.ga}}{\le}
\liminf_{n\to \infty} F_n(u_n, h_n - D_tu_n) 
\\
&\overset{\eqref{eq.par.min1}_2}{=}
\!\! \liminf_{n\to \infty} \; \int_0^T \langle h_n - D_tu_n, u_n\rangle \, d\mu(t) 
\overset{\eqref{eq.par.co1}}{=}
\!\int_0^T \! \langle h - D_tu, u\rangle \, d\mu(t) .  
\end{split}
\end{equation}
Thus, defining $\widetilde\Phi$ as in \eqref{eq.par.BENfun''} by simply dropping the index $n$,
\begin{equation}\label{eq.par.min2}  
\left\{\begin{split} 
&u\in {\cal V}, 
\\[1mm]
&\widetilde\Phi(u,h) = 0\;\; 
\big(\! =\inf_{\cal V} \widetilde\Phi(\cdot,h) \big),
\\
&u(0) = u^0.
\end{split}\right.
\end{equation} 
By Theorem~\ref{teo.repr'}, \eqref{pb=} is thus established. 
\hfill$\Box$ 

\begin{remarks}\rm
(i) This theorem improves the results of Section~8 of \cite{ViCalVar},
in particular by formulating the limit functional in terms of the memoryless function $\varphi$.

(ii) The reader will notice that we did not directly prove $\Gamma$-compactness and 
$\Gamma$-convergence of the BEN-type functionals $\widetilde\Phi_n$s. 
We rather derived existence of a solution of the initial-value problem, 
and then passed to the limit in the null-minimization problem.
(Loosely speaking, we used this procedure because of a mismatch between the space 
for which we get compactness and that for which we have stability. 
As one might expect, the term $\langle D_tv,v\rangle$ is the main source of difficulty. 
But it does not seem worth to illustrate this issue here.)

(iii) The question whether the representable operator $\widehat\alpha$ is maximal monotone 
is here left open. 
  
(iv) As we saw in Section~\ref{sec.BEN}, when dealing with a time-dependent maximal monotone operator $\alpha: V\!\times {}]0,T[{} \to {\cal P}(V')$, 
one may represent this operator by a time-dependent representative function 
$f_\alpha(\cdot,\cdot,t)\in {\cal F}(V)$, and then reformulate the initial-value problem
in terms of an extended BEN principle; see Theorem~\ref{teo.repr}.

Theorem~\ref{eq.par.stru} takes over to this more general set-up,
as an existence result analogous to Lemma~\ref{regul} holds for time-dependent 
maximal monotone operators.
In this case $\psi$ explicitly depends on time;
so the generality of Theorems~\ref{teo.comp} and \ref{teo.comp'} is fully exploited.
We leave the details to the interested reader.

(v) If instead of prescribing an initial condition one assumes the solution $u$
to be $T$-periodic in time, then $\int_0^T \langle D_tu,u\rangle \, dt =0$.
In this case it is not needed to introduce any weight function,
and the above argument is much simplified.
\end{remarks}

\bigskip
\noindent{\bf Applications to PDEs.\/}
If $\Omega$ is a bounded Lipschitz domain of $\erre^N$ ($N\ge 1$)
and $\vec\phi$ is a maximal monotone mappings $\erre^N\to {\cal P}(\erre^N)$,
one may take
\begin{equation}\label{eq.spaces.1}
H= L^2(\Omega), \quad
V= H^1_0(\Omega), \quad
\alpha(v) = - \nabla \!\cdot\! \vec\phi(\nabla v)  
\;\;\hbox{ in }{\cal D}'(\Omega). 
\end{equation}
If $N=3$, denoting the outward-oriented unit normal vector-field on
$\partial\Omega$ by $\vec\nu$, one may deal with
\begin{equation}\label{eq.spaces.2}
\begin{split}
&H= \big\{\vec v\in L^2(\Omega)^3: \nabla \!\cdot\! \vec v =0 \hbox{ a.e.\ in }\Omega)
\big\},
\\
&V= \big\{\vec v\in H: \nabla \!\times\! \vec v \in L^2(\Omega)^3, \;
\vec\nu \!\times\! \vec v = \vec 0 \hbox{ in } H^{-1/2}(\partial\Omega)^3 \big\},  
\\
&\vec\alpha(\vec v) = \nabla \times \vec\phi(\nabla \!\times\! \vec v) 
\quad\hbox{ in }{\cal D}'(\Omega)^3, \forall \vec v\in V.
\end{split}
\end{equation} 

One may also consider time-dependent operators.
In any of these cases one may apply Theorem~\ref{eq.par.stru} to the 
corresponding PDE
\begin{equation}  
D_tu + \alpha(u) \ni h
\qquad\hbox{ in }\Omega \!\times {}]0,T[,
\end{equation} 
or, after some modifications, e.g.\ to 
\begin{equation}  
D_tu - \Delta \int_0^t u(\cdot,\tau) \, d\tau + \alpha(u) \ni h
\qquad\hbox{ in }\Omega \!\times {}]0,T[.
\end{equation}

\section{Structural stability of doubly-nonlinear flows. I} 
\label{sec.DNE}   

\noindent 
In this section we study the structural compactness and stability of doubly-nonlinear equations of the form 
\begin{equation}\label{eq.DNE.1}
D_t\partial\gamma(u) + \alpha(u)\ni h, 
\end{equation} 
with $\alpha$ a maximal monotone operator and $\gamma$ a lower semicontinuous convex function.

Let $V$ and $H$ be real separable Hilbert spaces as in \eqref{eq.par.tripl}, and set 
\begin{equation}\label{eq.DNE.1+}
{\cal W} := L^\infty(0,T;H)\cap H^1(0,T;V').
\end{equation} 
Let four sequences $\{\alpha_n\}, \{\gamma_n\}, \{u^0_n\}$ and $\{h_n\}$ 
be given such that
\begin{eqnarray}
&\forall n, \; \alpha_n: V\to {\cal P}(V') \hbox{ is maximal monotone, } 
\label{eq.DNE.a1}
\\
&\exists C_1, C_2>0: \forall n, 
\forall (v,v^*)\in \graph(\alpha_n), \quad 
\langle v^*,v \rangle \ge C_1|v\|_V^2 - C_2, 
\label{eq.DNE.a2}
\\ 
&\exists C_3, C_4>0: \forall n, 
\forall (v,v^*)\in \graph(\alpha_n), \quad 
\|v^*\|_{V'} \le C_3 \|v\|_V + C_4,  
\label{eq.DNE.a3}
\\
&\alpha_n(0) \ni 0 \qquad\forall n,
\label{eq.DNE.a4}
\\
&\forall n, \; \gamma_n: H\to \erre \hbox{ is convex and lower semicontinuous, } 
\label{eq.DNE.a5}
\\
&\exists \bar C_1,..., \bar C_4>0: \forall n, \forall v\in H, \quad 
\bar C_1|v\|_H^2 - \bar C_2\le \gamma_n(v) \le \bar C_3 \|v\|_H^2 + \bar C_4,  
\label{eq.DNE.a6}
\\ 
&w^0_n\to w^0 \quad\hbox{ in }H,
\label{eq.DNE.a7}
\\
&h_n\to h \quad\hbox{ in }L^2(0,T;V').
\label{eq.DNE.a8}
\end{eqnarray} 
For any $n$, we formulate the following initial-value problem
\begin{equation}\label{eq.DNE.Cau1}  
\left\{\begin{split} 
&u_n\in L^2(0,T;V), \quad w_n\in {\cal W},
\\
&D_tw_n + \alpha_n(u_n)\ni h_n \quad \hbox{ in $V'$, a.e.\ in }{}]0,T[,
\\
&w_n\in \partial\gamma_n(u_n) \quad \hbox{ in $H$, a.e.\ in }{}]0,T[,
\\
&w_n(0) = w^0_n.
\end{split}\right.
\end{equation}  
 
\begin{lemma}[\cite{DiSh}, \cite{ViDNE}]\label{pro.DNE.exi}  
Under the hypotheses \eqref{eq.par.tripl}, \eqref{eq.DNE.a1}--\eqref{eq.DNE.a8}, 
for any $n$ the initial-value problem~\eqref{eq.DNE.Cau1} has at least one solution 
$(u_n,w_n)\in L^2(0,T;V) \!\times\! {\cal W}$.
Moreover, the sequence $\{(u_n,w_n)\}$ is bounded in this space.
\end{lemma}

\medskip
\noindent{\bf Null-Minimization.}
Next we reformulate the problem \eqref{eq.DNE.Cau1}. 
Let us first notice that 
\begin{equation}\label{eq.DEN.double} 
\begin{split}
\int_0^T \! d\tau\! \int_0^\tau \! \langle D_tw,v\rangle \, dt 
= \int_0^T \! d\tau\! \int_0^\tau D_t\gamma_n^*(w) dt 
= \int_0^T \gamma_n^*(w(\tau))\, d\tau - T\gamma_n^*(w(0)) &
\\
\qquad\forall v\in L^2(0,T;V),\forall w\in {\cal W},
\hbox{ with $\forall w\in \partial\gamma_n(v)$ a.e.\ in }]0,T[, &
\end{split}
\end{equation}
and define the measure $\mu$ as in \eqref{eq.comstab.mu'}.
For any $n$, let us represent the operator $\alpha_n$ by
$\varphi_n = (\pi+ I_{\alpha_n})^{**} \; (\in{\cal F}(V))$, define the affine subspace
\begin{equation}\label{eq.DNE.aa}
{\cal W}_{\mu,u_n^0} := \big\{v\in L^2(0,T;H): 
D_tv \in L^2_\mu(0,T;V'), v(0) = w_n^0\big\}, 
\end{equation}
and the nonnegative twice-time-integrated functional 
\begin{equation}\label{eq.DNE.fun1}
\begin{split}
&\widetilde\Phi_n(u,w,u^*) 
:=  \int_0^T \! d\tau\! \int_0^\tau \big[\gamma_n(u) + \gamma_n^*(w) 
-(w,u)_H \big] \, dt 
\\
&+ \bigg(\!\int_0^T \! d\tau\! \int_0^\tau \big[ \varphi_n(u,u^* -D_tw) 
- \langle u^*-D_tw,u\rangle \big] \, dt \bigg)^+
\\
&\overset{\eqref{eq.DEN.double}}{=}
\int_0^T \! \big[\gamma_n(u) + \gamma_n^*(w) -(w,u)_H \big] \, d\mu(t)
\\
&+ \bigg(\!\! \int_0^T \big[\varphi_n(u,u^* -D_tw) -\langle u^*,u\rangle \big] \, d\mu(t)  
+ \int_0^T \! \gamma_n^*(w(\tau))\, d\tau - T\gamma_n^*(w_n^0) \!\bigg)^+
\\
& \qquad \qquad \qquad\qquad\qquad\quad 
\forall u\in L^2_\mu(0,T;V), \forall w\in {\cal W}_{\mu,u_n^0}, 
\forall u^*\in L^2_\mu(0,T;V'),
\\[1mm]
& \widetilde\Phi_n(u,w,u^*)  
:= +\infty \quad\; \hbox{ for any other }(u,w,u^*) \in 
L^2(0,T;V) \!\times\! {\cal W} \!\times\! L^2(0,T;V').
\end{split}
\end{equation}
 
\begin{proposition}\label{pro.DNE.equivn} 
For any $n$, the pair $(u_n,w_n)$ solves the initial-value problem~\eqref{eq.DNE.Cau1} 
if and only if  
\begin{equation}\label{eq.DNE.min1}  
\left\{\begin{split} 
&u_n\in L^2(0,T;V), \quad w_n\in {\cal W},
\\
&\widetilde\Phi_n(u_n,w_n,h_n) = 0\;\; 
\big(\! =\inf_{L^2(0,T;V)\times {\cal W}} \widetilde\Phi_n(\cdot,\cdot,h_n) \big),
\\
&w_n(0) = w^0_n.
\end{split}\right.
\end{equation} 
\indent
Moreover, \eqref{eq.DNE.min1}$_2$ and \eqref{eq.DNE.min1}$_3$ are equivalent to 
\begin{equation}\label{eq.DNE.min1'}
\left\{\begin{split}
&\dps \int_0^T \big[\gamma_n(u_n) + \gamma_n^*(w_n) - (w_n, u_n)_H\big] \, d\mu(t)
\le 0,
\\
&\dps \int_0^T \big[\varphi_n(u_n,u_n^* -D_tw_n) -\langle u_n^*,u_n\rangle\big] \, d\mu(t)  
+\! \int_0^T \!\! \gamma_n^*(w_n(\tau))\, d\tau - T\gamma_n^*(w_n^0) \le 0.
\end{split}\right.
\end{equation} 
\end{proposition}

\medskip
\noindent{\bf Proof.\/} 
By the Fenchel system \eqref{eq.fitzp.Fitztheo2}, the first integrand of \eqref{eq.DNE.fun1} 
is nonnegative.
The null-minimization principle \eqref{eq.DNE.min1} is thus equivalent to the system
\eqref{eq.DNE.min1'}. 

The first inequality of \eqref{eq.DNE.min1'} is equivalent to \eqref{eq.DNE.Cau1}$_3$.
By adapting the procedure of the time-integrated BEN principle, see Theorem~\ref{teo.repr'}, 
and noting that by \eqref{eq.DNE.Cau1}$_3$ 
$\langle D_tw_n u_n\rangle = D_t \gamma_n^*(w_n)$,
we see that the second inequality of \eqref{eq.DNE.min1'} is tantamount to
\[
\int_0^T \! d\tau\! \int_0^\tau \big[ \varphi_n(u_n,u_n^* -D_tw_n) 
- \langle u_n^* -D_tw_n,u_n\rangle \big] \, dt \le 0.
\]
As $\varphi_n$ is a representative function, this is equivalent to \eqref{eq.DNE.Cau1}$_2$.
\hfill$\Box$ 

\begin{theorem}[Structural compactness and stability]\label{eq.DNE.stru}
Let \eqref{eq.par.tripl}, \eqref{eq.DNE.a1}--\eqref{eq.DNE.a8} be fulfilled.
For any $n$, let $(u_n,w_n)$ be a solution of problem~\eqref{eq.DNE.Cau1},
and define $\varphi_n = (\pi+ I_{\alpha_n})^{**} \; (\in{\cal F}(V))$. Then:
\\
\indent
(i) There exist $u\in L^2(0,T;V)$ and $w\in {\cal W}$ such that, 
possibly extracting a subsequence,  
\begin{eqnarray} 
&u_n\wto u \qquad\hbox{ in }L^2(0,T;V),
\label{eq.DNE.ucon}
\\
&u_n\wsto u \qquad\hbox{ in }{\cal W}.
\label{eq.DNE.wcon}  
\end{eqnarray} 
\indent
(ii) There exists a function $\varphi\in {\cal F}(V)$ such that, 
defining the measure $\mu$ as in \eqref{eq.comstab.mu'} and setting 
\begin{equation}\label{eq.DNE.de''}  
\begin{split} 
&\psi_{n,(v,v^*)}(t) = \varphi_n(v(t), v^*(t)), \qquad
\psi_{(v,v^*)}(t) = \varphi(v(t), v^*(t)),  
\\[1mm]
&\hbox{for a.e.\ }t\in{}]0,T[,\forall (v, v^*)\in L^2_\mu(0,T;V \!\times\! V'), \forall n,
\end{split} 
\end{equation}
then $\psi_n,\psi: L^2_\mu(0,T;V \!\times\! V') \to L^1_\mu(0,T)$ and, 
possibly extracting a subsequence, 
\begin{eqnarray}  
&\begin{split}
&\hbox{ $\psi_n$ sequentially $\Gamma$-converges to }\psi
\\
&\hbox{ in the topology $\widetilde\pi$ of $L^2_\mu(0,T;V \!\times\! V')$ and } 
\\
&\hbox{ in the weak topology of $L^1_\mu(0,T)$ (cf.\ \eqref{eq.evol.defgamma.2}). }
\end{split}
\label{eq.DNE.tesi+}
\end{eqnarray}
\indent
(iii) There exists a function $\gamma: H\to $ that fulfills lower and upper estimates analogous 
to \eqref{eq.DNE.a6}, and such that, possibly extracting a subsequence, 
\begin{eqnarray}  
&\begin{split}
&\hbox{ $\gamma_n$ strongly $\Gamma$-converges to $\gamma$ in $L^2(0,T;H)$, and }
\label{eq.DNE.r}
\\
&\hbox{ $\gamma_n^*$ sequentially weakly $\Gamma$-converges to $\gamma^*$ in
$L^2(0,T;H)$. }
\label{eq.DNE.s}
\end{split} 
\end{eqnarray}
\\
\indent
(iv) Denoting by $\alpha: V\to {\cal P}(V')$ 
the monotone operator that is represented by $\varphi$,
the pair $(u,w)$ solves the initial-value problem
\begin{equation}\label{eq.DNE.Cau=}  
\left\{\begin{split} 
&u\in L^2(0,T;V), \quad w\in {\cal W},
\\
&D_tw + \alpha(u)\ni h \quad \hbox{ in $V'$, a.e.\ in }{}]0,T[,
\\
&w\in \partial\gamma(u) \quad \hbox{ in $H$, a.e.\ in }{}]0,T[,
\\
&w(0) = w^0.
\end{split}\right.
\end{equation} 
\end{theorem}

\medskip
\noindent{\bf Proof.\/}
(i) By \eqref{eq.DNE.a5} and \eqref{eq.DNE.a6}, there exists
a convex and lower semicontinuous function $\gamma: V\to \erre$ such that
\begin{eqnarray}
&\gamma_n\Gammato \gamma  \quad\hbox{ strongly in }L^2(0,T;H),
\label{eq.DNE.a10}
\\
&\bar C_1|v\|_H^2 - \bar C_2\le \gamma(v) \le \bar C_3 \|v\|_H^2 + \bar C_4
\qquad\forall v\in H.  
\label{eq.DNE.a11} 
\end{eqnarray} 
After e.g.\ \cite{At} p.\ 282-283, this entails that 
\begin{equation}\label{eq.DNE.a12} 
\gamma_n^*\Gammato \gamma^* \quad\hbox{ weakly in }L^2(0,T;H).
\end{equation} 

(ii) By Lemma~\ref{pro.DNE.exi}, \eqref{eq.DNE.ucon} and \eqref{eq.DNE.wcon} hold up to extracting subsequences. By \eqref{eq.comstab.motiv3} then
\begin{equation}\label{eq.DNE.co1}
\int_0^T \langle u_n, h_n - D_tw_n\rangle \, d\mu(t) \to
\int_0^T \langle u,h - D_tw\rangle \, d\mu(t),
\end{equation}
so that
\begin{equation}\label{eq.DNE.co2}
(u_n,h_n - D_tw_n)\pito (u,h - D_tw)
\qquad\hbox{ in }L^2_\mu(0,T;V \!\times\! V').
\end{equation} 

(iii) Because of \eqref{eq.DNE.a2} and \eqref{eq.DNE.a3},
we may apply Theorem~\ref{teo.comp'}.
As here the functions $\varphi_n$'s are translation-invariant, 
thus there exists a function $\varphi\in {\cal F}(V)$ such that,
defining $\psi$ as in \eqref{eq.DNE.de''} and possibly extracting a subsequence, 
\eqref{eq.DNE.tesi+} is fuflilled. This yields
\begin{equation}\label{eq.DNE.ga} 
\begin{split}
F_n(v,v^*) := \int_0^T \! \varphi_n(v,v^*) \, d\mu(t) 
\Gammatopi \int_0^T \! \varphi(v,v^*) \, d\mu(t) =: F(v,v^*)
\\
\qquad\hbox{ sequentially in }L^2_\mu(0,T;V \!\times\! V').
\end{split}
\end{equation}

Therefore
\begin{equation}
\begin{split}
&F(u, h- D_tw) 
\overset{\eqref{eq.DNE.co2},\eqref{eq.DNE.ga}}{\le}
\liminf_{n\to \infty} F_n(u_n, h_n - D_tw_n) 
\\
&\overset{\eqref{eq.DEN.double},\eqref{eq.DNE.min1'}_2}{\le}
\!\! \liminf_{n\to \infty} \; \int_0^T \langle u_n, h_n - D_tw_n\rangle \, d\mu(t)
\overset{\eqref{eq.DNE.co2}}{=}
\int_0^T \langle u,h - D_tw\rangle \, d\mu(t).  
\end{split}
\end{equation}
Thus $u$ fulfills the twice-time-integrated BEN principle, i.e.\ \eqref{eq.DNE.Cau1}
here written without the index $n$. \eqref{eq.DNE.Cau=} is thus established. 
\hfill$\Box$

\begin{remarks}\rm
(i) The above results may be applied to several doubly-nonlinear parabolic PDEs.
For instance, let $\Omega$ be a bounded Lipschitz domain of $\erre^N$ ($N\ge 1$),
let $H,V,\alpha$ be as in \eqref{eq.spaces.1},  
and let $\gamma: \erre\to \erre\cup \{+\infty\}$ be a lower semicontinuous convex function.
One may apply Theorem~\ref{eq.DNE.stru} to the quasilinear equation
\begin{equation}  
D_t\partial\gamma(u) + \alpha(u) \ni h
\qquad\hbox{ in }\Omega \!\times {}]0,T[.
\end{equation}  
 
If $\gamma: \erre^3\to \erre\cup \{+\infty\}$ then one may also deal with the 
analogous equation in the framework of \eqref{eq.spaces.2}.

(ii) If instead of prescribing an initial condition one assumes $w$ 
to be $T$-periodic in time, then $\int_0^T \langle D_tw,u\rangle \, dt =0$
(as $u\in \partial\gamma^*(w)$).
In this case it is not needed to introduce any weight function, 
and the above argument is much simplified.
\end{remarks}

\section{Structural stability of doubly-nonlinear flows. II} 
\label{sec.DNE'}  
 
\noindent 
In this section we study the structural stability of doubly-nonlinear equations of the form 
\begin{equation}\label{eq.DNE'.1}
\alpha(D_tu) + \partial\gamma(u)\ni u^*,
\end{equation}
with $\alpha$ a maximal monotone operator and $\gamma$ a lower semicontinuous convex function.
The reader will notice analogies and also differences between these developments and those of Section~\ref{sec.DNE}.
 
Let the Hilbert spaces $V, H$ be defined as in \eqref{eq.par.tripl}, and set
\begin{equation}\label{eq.DNE'.1=} 
{\cal U} = H^1(0,T;H)\cap L^\infty(0,T;V).
\end{equation}
Let four sequences $\{\alpha_n\}, \{\gamma_n\}, \{u^0_n\}$ and $\{h_n\}$ 
be given such that
\begin{eqnarray}
&\forall n, \; \alpha_n: H\to {\cal P}(H) \hbox{ is maximal monotone, } 
\label{eq.DNE'.a1}
\\
&\exists C_1, C_2>0: \forall n, 
\forall (v,v^*)\in \graph(\alpha_n), \quad 
\langle v^*,v \rangle \ge C_1|v\|_H^2 - C_2, 
\label{eq.DNE'.a2}
\\ 
&\exists C_3, C_4>0: \forall n, 
\forall (v,v^*)\in \graph(\alpha_n), \quad 
\|v^*\|_H \le C_3 \|v\|_H + C_4,  
\label{eq.DNE'.a3}
\\
&\alpha_n(0) \ni 0 \qquad\forall n,
\label{eq.DNE'.a4}
\\
&\forall n, \; \gamma_n: V\to \erre \hbox{ is convex and lower semicontinuous, } 
\label{eq.DNE'.a5}
\\
&\exists \bar C_1,..., \bar C_4>0: \forall n, \forall v\in V, \quad 
\bar C_1|v\|_V^2 - \bar C_2\le \gamma_n(v) \le \bar C_3 \|v\|_V + \bar C_4,  
\label{eq.DNE'.a6}
\\
&u^0_n\to u^0 \quad\hbox{ in }V,
\label{eq.DNE'.a7}
\\
&h_n\to h \quad\hbox{ in }L^2(0,T;H).
\label{eq.DNE'.a8}
\end{eqnarray} 
For any $n$, we formulate the following initial-value problem
\begin{equation}\label{eq.DNE'.Cau1}  
\left\{\begin{split} 
&u_n\in {\cal U}, \quad z_n\in L^\infty(0,T;H),
\\
&\alpha_n(D_tu_n) + z_n\ni h_n \quad \hbox{ in $H$, a.e.\ in }{}]0,T[,
\\
&z_n\in \partial\gamma_n(u_n)  \quad \hbox{ in $V'$, a.e.\ in }{}]0,T[,
\\
&u_n(0) = u^0_n.
\end{split}\right.
\end{equation} 
 
\begin{lemma}[\cite{CoVi}, \cite{ViDNE}] \label{pro.DNE'.exi}  
Under the hypotheses \eqref{eq.par.tripl}, \eqref{eq.DNE'.a1}--\eqref{eq.DNE'.a8}, 
for any $n$ the initial-value problem~\eqref{eq.DNE'.Cau1} has at least one solution 
$(u_n,z_n)\in {\cal U} \!\times\! L^\infty(0,T;H) $.
Moreover, the sequence $\{(u_n,z_n)\}$ is bounded in this space.
\end{lemma}

\medskip
\noindent{\bf Null-Minimization.}
Next we reformulate the problems~\eqref{eq.DNE'.Cau1}. 
Let us first notice that 
\begin{equation}\label{eq.DEN'.double} 
\begin{split}
\int_0^T \! d\tau\! \int_0^\tau \! (D_tv,z)_H \, dt 
= \!\int_0^T \!\! d\tau\!\! \int_0^\tau \! D_t\gamma_n(v(\tau)) dt
= \!\int_0^T \!\! \gamma_n(v(\tau))\, d\tau - T\gamma_n(v(0)) &
\\
\qquad\forall v\in {\cal U},\forall w\in {\cal W},
\hbox{ with $\forall z\in \partial\gamma_n(v)$ a.e.\ in }]0,T[. &
\end{split}
\end{equation}
For any $n$, let us represent the operator $\alpha_n$ by
$\varphi_n = (\pi+ I_{\alpha_n})^{**} \; (\in{\cal F}(V))$,
define the measure $\mu$ as in \eqref{eq.comstab.mu'}, the affine subspace
\begin{equation}\label{eq.DNE'.sp}
{\cal U}_{\mu,u_n^0} := \big\{v\in {\cal U}: v(0) = u_n^0\big\},
\end{equation}
and the nonnegative twice-time-integrated functional 
\begin{equation}\label{eq.DNE'.fun1}
\begin{split}
&\widetilde\Phi_n(u,z,u^*) 
:=  \int_0^T \! d\tau\! \int_0^\tau \big[\gamma_n(u) + \gamma_n^*(z) 
- \langle z,u\rangle \big] \, dt 
\\
&+ \bigg(\!\int_0^T \! d\tau\! \int_0^\tau \big[\varphi_n(D_tu,u^* -z)
- (D_tu,u^* -z)_H \big] \, dt \bigg)^+
\\
&\overset{\eqref{eq.DEN'.double}}{=}
\int_0^T \! \big[\gamma_n(u) + \gamma_n^*(z) - \langle z,u\rangle \big] \, d\mu(t)
\\
&+ \!\bigg(\!\int_0^T \! \big[\varphi_n(D_tu,u^* -z) - (D_tu,u^*)_H \big] \, d\mu(t) 
+ \!\!\int_0^T \!\!\! \gamma_n(u(\tau))\, d\tau - T\gamma_n(u(0)) \!\bigg)^+
\\
& \qquad \qquad  \qquad \qquad\qquad 
\forall u\in {\cal U}_{\mu,u_n^0}, \forall z\in L^2_\mu(0,T;H),
\forall u^*\in L^2_\mu(0,T;V'),
\\[1mm]
& \widetilde\Phi_n(u,z,u^*)  
:= +\infty \qquad  \hbox{ for any other }(u,z,u^*) \in 
{\cal U} \!\times\! {\cal W} \!\times\! {\cal V}'.
\end{split}
\end{equation}

\begin{proposition}\label{pro.DNE'.equivn} 
For any $n$, the pair $(u_n,z_n)$ solves the initial-value problem~\eqref{eq.DNE'.Cau1} 
if and only if  
\begin{equation}\label{eq.DNE'.min1}  
\left\{\begin{split} 
& u_n\in {\cal U}, \quad z_n\in L^\infty(0,T;H),
\\
&\widetilde\Phi_n(u_n,z_n,h_n) = 0
\;\; \big(\! =\inf_{{\cal U} \!\times\! L^\infty(0,T;H)} 
\widetilde\Phi_n(\cdot,\cdot,h_n) \big),
\\
&u_n(0) = u^0_n.
\end{split}\right.
\end{equation} 
\indent 
Moreover, \eqref{eq.DNE'.min1}$_2$ and \eqref{eq.DNE'.min1}$_3$ are equivalent to 
\begin{equation}\label{eq.DNE'.min1'}
\left\{\begin{split}
&\dps \int_0^T \!\big[\gamma_n(u_n) + \gamma_n^*(z_n) - \langle z_n, u_n\rangle\big] 
\, d\mu(t) \le 0,
\\
&\dps \int_0^T \big[\varphi_n(D_tu_n, h_n -z_n) - (h_n, D_tu_n)_H\big] \, d\mu(t) 
+ \int_0^T \!\gamma(u_n(\tau)) \, d\tau - T\gamma(u_n^0) \le 0.
\end{split}\right.
\end{equation} 
\end{proposition}

\medskip
\noindent{\bf Proof.\/} 
By the Fenchel system \eqref{eq.fitzp.Fitztheo2}, the first integrand of \eqref{eq.DNE'.fun1} 
is nonnegative.
The null-minimization principle \eqref{eq.DNE'.min1} is thus equivalent to the system
\eqref{eq.DNE'.min1'}. 

The first inequality of \eqref{eq.DNE'.min1'} is equivalent to \eqref{eq.DNE'.Cau1}$_3$.
This entails that 
\[
\int_0^T (D_tu_n,z_n)_H \, dt =  \gamma(u_n(T)) - \gamma(u_n^0),
\]
and this yields \eqref{eq.DNE'.Cau1}$_4$.
The second inequality of \eqref{eq.DNE'.min1'} is then equivalent to 
\[
\int_0^T \big[\varphi_n(D_tu_n, h_n -z_n) - (D_tu_n, h_n -z_n)_H \big] \, d\mu(t) \le 0,
\]
which is tantamount to \eqref{eq.DNE'.Cau1}$_2$.
\hfill$\Box$

\begin{theorem}[Structural compactness and stability]\label{eq.DNE'.stru}
Let \eqref{eq.par.tripl}, \eqref{eq.DNE'.a1}--\eqref{eq.DNE'.a8} be fulfilled.
For any $n$, let $(u_n,z_n)$ be a solution of problem~\eqref{eq.DNE'.Cau1},
and set $\varphi_n = (\pi+ I_{\alpha_n})^{**} \; (\in{\cal F}(V))$. Then: 
\\
\indent
(i) There exist $u\in {\cal U}$ and $z\in L^\infty(0,T;H)$ such that, 
possibly extracting a subsequence,  
\begin{eqnarray} 
&u_n\wsto u \qquad\hbox{ in }{\cal U},
\label{eq.DNE'.ucon}
\\
&z_n\wsto z \qquad\hbox{ in }L^\infty(0,T;H).
\label{eq.DNE'.wcon}  
\end{eqnarray} 
\indent
(ii) There exists a function $\varphi\in {\cal F}(H)$ such that, setting 
\begin{equation}\label{eq.DNE'.de''} 
\begin{split} 
&\psi_{n,(v,v^*)}(t) = \varphi_n(v(t),v^*(t)), \qquad
\psi_{(v,v^*)}(t) = \varphi(v(t),v^*(t)),  
\\[1mm]
&\hbox{for a.e.\ }t\in{}]0,T[,\forall (v,v^*)\in L^2_\mu(0,T;H^2), \forall n,
\end{split} 
\end{equation}
then $\psi_n,\psi: L^2_\mu(0,T;H \!\times\! H) \to L^1_\mu(0,T)$ and, 
possibly extracting a subsequence, 
\begin{eqnarray}  
&\begin{split}
&\hbox{ $\psi_n$ sequentially $\Gamma$-converges to }\psi
\\
&\hbox{ in the topology $\widetilde\pi$ of $L^2_\mu(0,T;H \!\times\! H)$ and } 
\\
&\hbox{ in the weak topology of $L^1_\mu(0,T)$ (cf.\ \eqref{eq.evol.defgamma.2}). }
\end{split}
\label{eq.DNE'.tesi+}
\end{eqnarray} 
\indent
(iii) There exists a function $\gamma: V\to \erre$ that fulfills lower and upper estimates 
analogous to \eqref{eq.DNE'.a6}, such that, possibly extracting a subsequence, 
\begin{eqnarray}  
&\begin{split}
&\hbox{ $\gamma_n$ strongly $\Gamma$-converges to $\gamma$ in $L^2(0,T;H)$, and }
\label{eq.DNE'.r}
\\
&\hbox{ $\gamma_n^*$ sequentially weakly $\Gamma$-converges to $\gamma^*$ in
$L^2(0,T;H)$. }
\label{eq.DNE'.s}
\end{split} 
\end{eqnarray}
\\
\indent
(iv) Denoting by $\alpha: H\to {\cal P}(H)$ 
the monotone operator that is represented by $\varphi$,
the pair $(u,z)$ solves the corresponding initial-value problem
\begin{equation}\label{eq.DNE'.Cau1=}
\left\{\begin{split} 
&u\in {\cal U}, \quad z\in L^\infty(0,T;H),
\\
&\alpha(D_tu) + z\ni h \quad \hbox{ in $H$, a.e.\ in }{}]0,T[,
\\
&z\in \partial\gamma(u)  \quad \hbox{ in $V'$, a.e.\ in }{}]0,T[,
\\
&u(0) = u^0.
\end{split}\right.
\end{equation}
\end{theorem}

\medskip
\noindent{\bf Proof.\/} 
(i) Because of \eqref{eq.DNE'.a2} and \eqref{eq.DNE'.a3},
we may apply Theorem~\ref{teo.comp'}.
As here the functions $\varphi_n$'s are translation-invariant, 
there exists a function $\varphi\in {\cal F}(V)$ such that,
defining $\psi$ as in \eqref{eq.DNE'.de''} and possibly extracting a subsequence, 
\eqref{eq.DNE'.tesi+} is fulfilled. 

By \eqref{eq.DNE'.a5} and \eqref{eq.DNE'.a6}, there exists
a convex and lower semicontinuous function $\gamma: V\to \erre$ such that
\begin{eqnarray}
&\gamma_n\Gammato \gamma  \qquad\hbox{ strongly in }L^2(0,T;H),
\label{eq.DNE'.a10}
\\
&\bar C_1|v\|_H^2 - \bar C_2\le \gamma(v) \le \bar C_3 \|v\|_H^2 + \bar C_4
\qquad\forall v\in H.
\label{eq.DNE'.a11} 
\end{eqnarray} 
After e.g.\ \cite{At} p.\ 282-283, this entails that 
\begin{equation}\label{eq.DNE'.a12} 
\gamma_n^*\Gammato \gamma^* \qquad\hbox{ weakly in }L^2(0,T;H).
\end{equation} 

(ii) By Lemma~\ref{pro.DNE'.exi}, \eqref{eq.DNE'.ucon} and \eqref{eq.DNE'.wcon} 
hold up to extracting subsequences. This yields
\begin{equation}\label{eq.DNE'.s}  
u_n\to u \qquad\hbox{ in }L^2(0,T;H).
\end{equation} 

By \eqref{eq.DNE'.a10} and \eqref{eq.DNE'.a12}, 
passing to the inferior limit in \eqref{eq.DNE'.min1'}$_1$ we thus get
\begin{equation}
\int_0^T \big[\gamma(u) + \gamma^*(z)\big] \, dt
\le \int_0^T \langle z, u\rangle \, dt,
\end{equation} 
which is equivalent to
\begin{equation}\label{eq.DNE'.v}  
z\in \partial\gamma(u)  \quad \hbox{ in $H$, a.e.\ in }{}]0,T[.
\end{equation}
 
(iii) Notice that
\begin{equation}\label{eq.comstab.mot}
\begin{split}
&\int_0^T (D_tu_n,z_n)_H \, d\mu(t) 
\overset{\eqref{eq.DNE'.Cau1}_3}{=} \int_0^T D_t \gamma_n(u_n) \, (T-t) \, dt 
\\
&= \int_0^T \! \gamma_n(u_n) \, dt - T\gamma_n(u_n^0) 
\to \int_0^T \! \gamma(u) \, dt - T\gamma(u^0) 
\overset{\eqref{eq.DNE'.v}}{=}\int_0^T \! (D_tu,z)_H \, d\mu(t) .
\end{split}
\end{equation}
By \eqref{eq.DNE'.a8} then 
\begin{equation}\label{eq.DNE'.co}
(D_tu_n, h_n - z_n)\pito (D_tu, h - z)
\qquad\hbox{ in }L^2_\mu(0,T;H \!\times\! H).
\end{equation} 

(iv) \eqref{eq.DNE'.tesi+} yields
\begin{equation}\label{eq.DNE'.ga} 
\begin{split}
G_n(v,w) := \int_0^T \! \varphi_n(v,w) \, d\mu(t) 
\Gammatopi \int_0^T \! \varphi(v,w) \, d\mu(t) =: G(v,w)
\\
\qquad\hbox{ sequentially in }L^2_\mu(0,T;H \!\times\! H).
\end{split}
\end{equation}
 
Therefore
\begin{equation}
\begin{split}
&G(D_tu,h - z) 
\overset{\eqref{eq.DNE'.co},\eqref{eq.DNE'.ga}}{\le}
\!\! \liminf_{n\to \infty} G_n(D_tu_n, h_n - z_n)
\\
&\overset{\eqref{eq.DEN'.double},\eqref{eq.DNE'.min1'}_2}{\le} 
\!\! \liminf_{n\to \infty} \; \int_0^T \langle D_tu_n, h_n - z_n\rangle \, d\mu(t) 
\overset{\eqref{eq.DNE'.co}}{=} 
\int_0^T \langle D_tu,h - z\rangle \, d\mu(t).  
\end{split}
\end{equation}
Thus $u$ fulfills the time-integrated BEN principle, i.e.\ Theorem~\ref{teo.repr'}, 
\eqref{pb=} is thus established. 
\hfill$\Box$ 

\begin{remarks}\rm
(i) Theorem~\ref{eq.DNE'.stru} may be applied to several doubly-nonlinear parabolic PDEs.
For instance, let $\Omega$ be a bounded Lipschitz domain of $\erre^N$ ($N\ge 1$),
let the spaces $H,V,$ be as in \eqref{eq.spaces.1}, let
$\alpha: \erre^N\to {\cal P}(\erre^N)$ be maximal monotone,
and let $\gamma: \erre^N\to \erre\cup \{+\infty\}$ be a lower semicontinuous convex function.
One may apply Theorem~\ref{eq.DNE'.stru} to the quasilinear equation
\begin{equation}  
\alpha(D_tu) - \nabla \!\cdot\! \partial \gamma(\nabla u) \ni h
\qquad\hbox{ in }\Omega \!\times {}]0,T[.
\end{equation}  

We leave to the reader to identify the modifications that are needed to deal with an analogous equation in the set-up of \eqref{eq.spaces.2}. 

(ii) If instead of prescribing an initial condition one assumes $u$ 
to be $T$-periodic in time, then $\int_0^T \langle z,D_tu\rangle \, dt =0$
(as $z\in \partial\gamma(u)$).
In this case it is not needed to introduce any weight function, 
and the above argument is much simplified. 
\end{remarks}\rm

\bigskip\bigskip
\centerline{\bf Conclusions} 
\medskip
 
We introduced a definition of evolutionary $\Gamma$-convergence of weak type.
Via this notion and a nonlinear topology of weak type,
we proved the structural compactness and the structural stability of the weak formulation of the initial-value problem for several quasi-nonlinear PDEs.

The methods that we illustrated might be used to prove analogous structural properties of other evolutionary equations, including first-order flows for pseudo-monotone operators.

\bigskip\bigskip
\centerline{\bf Acknowledgment}
\medskip

The present research was partially supported by a MIUR-PRIN 10-11 grant for the project 
``Calculus of Variations" (Protocollo 2010A2TFX2-007).

\baselineskip=12truept
\bigskip

Author's address: 
\medskip

Augusto Visintin \par 
Universit\`a degli Studi di Trento \par 
Dipartimento di Matematica \par 
via Sommarive 14, \ 38050 Povo (Trento) - Italia \par 
Tel   +39-0461-281635 (office), +39-0461-281508 (secretary) \par 
Fax      +39-0461-281624 \par 
Email:   augusto.visintin@unitn.it 

\end{document}